\documentclass[a4paper]{article}
\usepackage[english]{babel}
\usepackage[margin=3cm]{geometry}
\usepackage{amssymb,amsmath,xcolor,graphicx,authblk}
\usepackage[hyphens]{url}

\graphicspath{{figures/}}

\newcommand{\R}{\mathbb{R}}

\newcommand{\A}{\mathcal{A}}

\newcommand{\tpose}{^\mathsf{T}}
\DeclareMathOperator{\Tr}{Tr}

\DeclareMathOperator*{\argmax}{argmax}

\title{Reliable optimal controls for SEIR models in epidemiology}
\date{26 April 2024}
\author{Simone Cacace}
\author{Alessio Oliviero}
\affil{Department of Mathematics, Sapienza University of Rome, Italy}

\begin{document}

\maketitle

\begin{abstract}
We present and compare two different optimal control approaches applied to SEIR models in epidemiology, which allow us to obtain some policies for controlling the spread of an epidemic. The first approach uses Dynamic Programming to characterise the value function of the problem as the solution of a partial differential equation, the Hamilton-Jacobi-Bellman equation, and derive the optimal policy in feedback form. The second is based on Pontryagin's maximum principle and directly gives open-loop controls, via the solution of an optimality system of ordinary differential equations. This method, however, may not converge to the optimal solution. We propose a combination of the two methods in order to obtain high-quality and reliable solutions. Several simulations are presented and discussed, also checking first and second order necessary optimality conditions for the corresponding numerical solutions.
\end{abstract}

\bigskip



\section{Introduction}
\label{sec:introduction}

The development of the Covid-19 pandemic has increased the interest towards the mathematical modelling of epidemics and, in the last three years, a number of papers dealing with several aspects of epidemiology have appeared. The number of topics is very large and ranges from the analytical description of infections to the socio-economic impact of the pandemic.  However, epidemiological models have been introduced even before the Covid-19 pandemic (see e.g. \cite{DH_2000} and the monography \cite{B_2011}) as a useful tool to analyse and predict the development of an infective disease. In particular, compartmental models describing the interactions between susceptible, infected and recovered consist of a system of ordinary differential equations, where the variables represent the categories which the population is divided into.  One of the simplest model, the popular SIR model, dates back to a paper by Kermack and McKendrick (1927) \cite{KMcK_1927}. Nowadays the modelling can include several aspects of an epidemic, such as space diffusion, age distribution, vaccination and local effects, and models are becoming more and more complex.

This work is focused on the control of epidemiological models. More precisely, given the description of the infection, we introduce some control parameters in the model, accounting for the possibility to get a vaccine, restrict social interactions or regulate the inflow of people from abroad. Clearly, by varying these parameters we will modify the evolution of the infection and its distribution in the population. In order to compare these effects, we introduce a cost functional associated with the development of the epidemic. This functional can contain different components, measuring costs in terms of victims, hospitalization procedures (such as intensive care) and economic impact. The final goal is to determine how restrictive measures and vaccines can affect the evolution of the disease and possibly find the optimal controls that minimise the cost at each time. To this end, we adopt two classical approaches based on general results in optimal control theory.
For the solution of an optimal control problem, one can look for controls that are expressed either as a function of time (open-loop controls) or as a function of the state of the system (feedback controls), depending on the method that is used to solve the problem. 

The first approach is based on Pontryagin's maximum principle (PMP), a very well-known variational framework which gives necessary conditions for a control to be optimal. However, the corresponding open-loop controls cannot react to perturbations of the system or to uncertainties in the model, since they are just a function of time, and this might be an issue in some cases. For a general introduction to this theory we refer to \cite{PBGM_1964,FR_1975,BCD_2008}.

The second approach follows the celebrated Dynamic Programming Principle (DPP) due to R. Bellman \cite{B_1957}, introducing the value function, a function of the state variables of the system, corresponding to the best cost among all the possible control strategies. This function can be characterised as the unique viscosity solution \cite{CL_1983,CEL_1984} of a non-linear partial differential equation, the Hamilton-Jacobi-Bellman (HJB) equation \cite{E_1998, BCD_2008}. Using a suitable numerical scheme of semi-Lagrangian type \cite{FG_1999}, one can calculate the value function and obtain the optimal controls. General results on this approach can be found in \cite{BCD_2008, FF_SIAM2013,FF_2016}. The main advantage of this strategy, despite the high computational cost, is that it directly gives feedback controls.

There is a huge amount of literature on optimal control applied to epidemiological models, aiming at describing different infectious diseases in different settings. These works, some of which date back to several years before the Covid-19 pandemic, are mainly based on Pontryagin's principle and on the solution of the associated optimality system \cite{Beh_2000, HD_2011, YB_2012, AFG_2022, FCAAHT_2023}. Moreover, the numerical schemes for the solution of the optimality system are typically based on \textit{local} descent methods, such as Newton iterations \cite{KS_1987, GDA_2020, IS_2013} or the so-called Direct-Adjoint Looping algorithm, also known as Forward-Backward Sweep \cite{XXXY_2017, BMO_2019, BDPV_2021}. To our knowledge, one of the less explored aspects of these optimisation techniques is that – if they converge – they may not converge to the global minimum of the cost functional \cite{DO_2021}, unless additional hypotheses are verified or the system is rather simple.

On the other hand, methods based on the Dynamic Programming approach are still scarcely employed for real-world problems, which involve complex or high-dimensional systems, due to their high computational cost. At present, the application of this approach to the control of epidemiological models is more theoretical or limited to low-dimensional settings \cite{L_1979, BGY_1997, GS_2014, LT_2015, STW_2019, SH_2019, HWL_2020, CLL_2023}.  However, the most relevant feature of this approach is the theoretical guarantee of convergence to the optimal solution, i.e. to the global minimum of the cost. In addition, the latest advancements in both CPU and GPU architectures make it possible to approximate the solution of large scale problems in a reasonable time, especially resorting to parallel computing techniques.

The present work aims at exploiting Dynamic Programming to improve and validate the results obtained with the variational approach. In particular, we combine the two approaches, using the optimal controls given by a semi-Lagrangian scheme as a warm guess to initialise the Direct-Adjoint Looping (DAL) method. Moreover, we evaluate the reliability of the optimal controls produced with the combined SL-DAL scheme by checking first and second order necessary optimality conditions along the computed numerical solutions. As it will be confirmed by our numerical experiments, especially when the considered epidemiological model gets more complicated, the DAL algorithm alone could provide various sub-optimal solutions, depending on the initial guess for the controls, but this does not happen if we initialise it with the output of the semi-Lagrangian scheme. Similar ideas, exploiting the theoretical connection between PMP and DPP, have been investigated in different settings: in \cite{CM_2010}, a rough finite difference approximation of the value function at single points is used for the initialization of a Newton-based shooting method for solving the PMP optimality system; in \cite{AKK_2021,YD_2021} PMP is interpreted as a representation formula for the solution of HJB equations, and it provides an efficient way to reduce space complexity when solving high dimensional control problems.

This paper is structured as follows. In Section \ref{sec:SEIR}, we present the general setting of the SEIR compartmental model, then we introduce the control parameters representing restrictive measures and vaccination, thus obtaining a \textit{controlled} SEIR model. In Section \ref{sec:cost}, we briefly outline the general theory of finite horizon optimal control problems and we formulate the cost functional. In Section \ref{sec:dp}, we introduce the Dynamic Programming approach, the related  Hamilton-Jacobi-Bellman equation and its semi-Lagrangian numerical approximation. Section \ref{sec:pontryagin} is dedicated to the variational approach to finite horizon problems, via Pontryagin's maximum principle. In particular, we describe the Direct-Adjoint Looping method to approximate the solution of the optimality system associated with the control problem, together with the corresponding first and second order necessary optimality conditions. Section \ref{sec:SLvsDAL} contains a more in depth comparison between the two approaches, highlighting their differences and their strengths, and also the proposed combination of the two methods. In Section \ref{sec:numerical}, we report the results of several numerical experiments performed with both methods on some variations of the controlled SEIR model. In particular, we consider the case of a control that can open or close the borders in a setting where there is a constant source term, representing incoming people from abroad. 
Finally, we also consider some state constraints that represent the limited availability of intensive care units, together with temporary immunity.

\section{The controlled SEIR model}
\label{sec:SEIR}

In this section, we present one of the most simple and frequently used compartmental models, the epidemic SEIR model, and an optimal control problem associated with it. Since our aim is to formulate an optimal control problem, we will not go through all the analytical properties of the model, and refer to \cite{AS_1984,DH_2000,DHB_2013} for further details.

\begin{figure}[h]
\centering
\includegraphics[width=9cm]{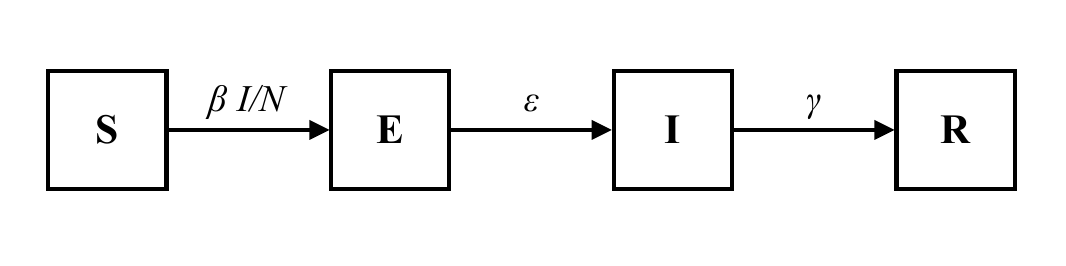}
\caption{Diagram of the epidemic SEIR model}
\label{fig:diagramma_seir}
\end{figure}

The population - which is assumed to be homogeneous and well mixed - is divided into four categories:
\begin{itemize}
    \item S: the susceptible portion, individuals who can contract the disease,
    \item E: the exposed portion, individuals who have been infected but are not yet infectious,
    \item I: the infective portion, infected individuals who can transmit the disease,
    \item R: the recovered portion, those who have acquired immunity.
\end{itemize}
An individual can belong to any of the categories, but only to one of them at each time. If we denote the total population with $N$, the model can be sketched as in Figure \ref{fig:diagramma_seir}.

The infection process is modeled with the usual law of mass action: the force of infection is $I/N$; $\beta$ is the transmission rate, namely the product of the average number of contacts per unit time for each individual and the probability of transmission per contact; $1/\varepsilon$ is the mean duration of the latency period; $1/\gamma$ is the mean infectious period. The corresponding differential model is therefore:
\begin{equation} \label{seir}
\begin{cases}
S' = - \beta S I/N \\
E' = \beta S I/N - \varepsilon E \\
I' = \varepsilon E - \gamma I \\
R' =  \gamma I
\end{cases}
\end{equation}
In order to avoid working with large numbers, we can normalise the equations by dividing everything by $N$, thus obtaining 
\begin{equation} \label{seir_norm}
\begin{cases}
s' = - \beta s i \\
e' = \beta s i - \varepsilon e \\
i' = \varepsilon e - \gamma i \\
r' =  \gamma i
\end{cases}
\end{equation}
with initial data $(s(0),e(0),i(0),r(0))=(s_0,e_0,i_0,r_0)$ such that $s_0+e_0+i_0+r_0=1$. Usually, as it is reasonable to assume, $r_0=0$ and $s_0 \approx 1$. In the normalised model, each variable represents the fraction of population within a certain compartment, instead of the absolute number of individuals. We point out that the right hand side of both \eqref{seir} and \eqref{seir_norm} sums to zero, implying the conservation of the total population. Indeed, setting $x:=(s,e,i,r)$, it readily follows that the set
$$
\mathcal{S}=\left\{ x \in [0,1]^4\ \big\lvert \ ||x||_1=1 \right\}
$$
is positively invariant for system \eqref{seir_norm}.

Finally, we observe that the fourth equation is independent from the previous three, hence we can reduce the system to
\begin{equation} \label{seir_ridotto}
\begin{cases}
s' = - \beta s i \\
e' = \beta s i - \varepsilon e \\
i' = \varepsilon e - \gamma i
\end{cases}
\end{equation}
and then exploit the conservation property to obtain $r(t)=1-s(t)-e(t)-i(t)$.

As in the simpler SIR model, the basic reproduction number for this model is 
$$
R_0= \frac{\beta}{\gamma}\,.
$$
This may seem surprising, but in fact, as long as we do not consider natural deaths, the duration of the latency period has no influence on $R_0$ \cite{DHB_2013}.

\begin{figure}[h]
    \centering
    \includegraphics[width=9cm]{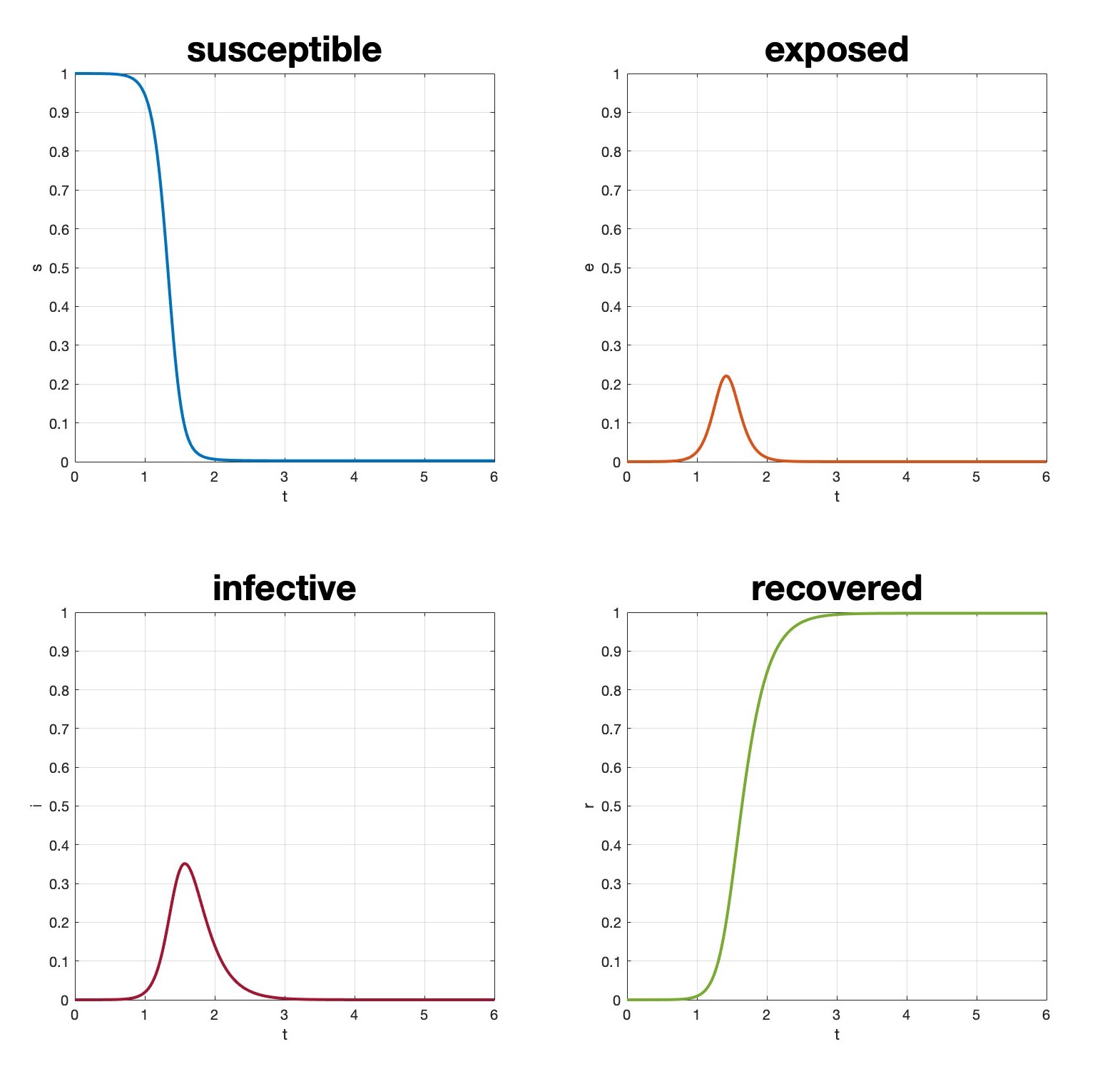}
    \caption{Trajectories of the normalised epidemic SEIR model. From left to right, top to bottom: susceptible, exposed, infective, recovered.}
    \label{fig:grafici_seir}
\end{figure}

What is peculiar of this model (see Figure \ref{fig:grafici_seir}) is the peak of infective, which is what the authorities usually want to minimise – or at least “spread” over a larger period of time – in order to avoid overloading health facilities.

\begin{figure}[h]
\centering
\includegraphics[width=9cm]{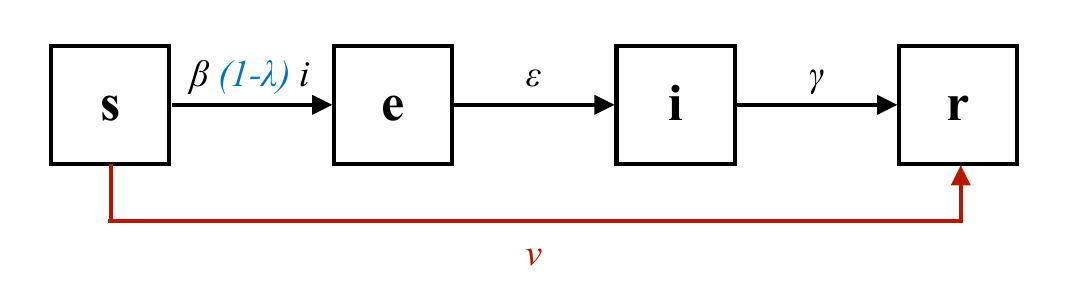}
\caption{Diagram of the normalised epidemic SEIR model with restrictions $\lambda$ and vaccinations $\nu$.}
\label{fig:diagramma_seir_controllato}
\end{figure}

Taking inspiration from the measures adopted to fight the SARS-CoV-2 pandemic, we introduce two types of controls that can act on the system at each time: restrictive measures ($\lambda$) and vaccinations ($\nu$). The two controls have quite different behaviours: while restrictive measures have the effect of reducing the contacts among the individuals, vaccines cut down the number of susceptible, making them immune without undergoing the infectious period. Considering directly the normalised model and assuming that only the susceptible get vaccinated, the resulting model can be sketched as shown in Figure \ref{fig:diagramma_seir_controllato}.
We multiply by $(1-\lambda)$ wherever the transmission rate $\beta$ appears in \eqref{seir_norm} and we add a new term containing $\nu$ which takes “mass” from the susceptible and moves it into the recovered, without passing through the other compartments. The model we obtain is thus the following:
\begin{equation} \label{seir_controllato}
\begin{cases}
s' = - \beta (1-\lambda) s i - \nu s \\
e' = \beta  (1-\lambda) s i - \varepsilon e \\
i' = \varepsilon e - \gamma i \\
r' =  \gamma i + \nu s
\end{cases}
\end{equation}
with the same initial data as before. In this way, the choice $\lambda=0$ corresponds to not modifying the contact rate at all, while $\lambda=1$ deletes the transmission term completely, meaning complete lock-down. Similarly, $\nu$ can be seen as a vaccination rate, where $\nu=0$ corresponds to not vaccinating anyone. We remark that there is no a priori upper bound $\nu_{\max}$ for the vaccination rate, as it depends on the resources available at each time in the society one is considering. Clearly, \eqref{seir_controllato} can be reduced as well and from now on we will always consider the following differential system:
\begin{equation} \label{seir_controllato_ridotto}
\begin{cases}
s' = - \beta (1-\lambda) s i - \nu s \\
e' = \beta  (1-\lambda) s i - \varepsilon e \\
i' = \varepsilon e - \gamma i
\end{cases}
\end{equation}
with initial data $(s(0),e(0),i(0))=(s_0,e_0,i_0)$ such that $s_0+e_0+i_0=1$. 

In Section \ref{sec:numerical} we will also consider some extensions of the basic SEIR model, taking into account state constraints, temporary immunity and incoming individuals from abroad.


\section{Finite horizon problems and cost definition}
\label{sec:cost}

In this section, we will briefly recall the theoretical setting for a \textit{finite horizon optimal control problem} and then define it in our particular case. To keep the notation lighter, we set $x:=(s_0,e_0,i_0) \in [0,1]^3$, $\alpha:=(\lambda,\nu) \in A$, where $A$ is a compact subset of $\R^2$, and consider a fixed time interval $[t,T]$. We will also refer to system \eqref{seir_controllato_ridotto} more concisely as
\begin{equation}\label{eq}
\left\{ \begin{array}{l}
\dot{y}(\tau)=f(y(\tau),\alpha(\tau),\tau), \;\; \tau\in(t,T],\\
y(t)=x,
\end{array} \right.
\end{equation}
where $y:[t,T] \to [0,1]^3$ represents the state variables and $f: [0,1]^3 \times A \times (t,T] \to [0,1]^3$ is the controlled dynamics. The expression “finite horizon” comes from the fact that $T \in \R$ is fixed, whereas for other kind of problems one can be interested in controlling the system as $T\to+\infty$. Since we are dealing with an epidemic model, the time scale is relatively short, months or at most a few years, therefore it is natural to fix a finite time horizon. The optimal control problem itself consists in finding a control function $\alpha^*$ in the space of \textit{admissible controls}
$$
\A := \left\{ \alpha: [t,T] \to A \ |\ \alpha \text{ measurable} \right\},
$$
minimising a given functional that represents the cost associated with the evolution of the epidemic. In general, the cost functional is of the form 
\begin{equation}\label{costo_orizz_fin}
J_{x,t}(\alpha) = \int_{t}^T \ell(y_{x,t} (\tau),\alpha(\tau),\tau)\ d\tau\ +\ g(y_{x,t}(T))
\end{equation}
where $y_{x,t} (\, \cdot \,)$ indicates the solution of \eqref{eq} starting at time $t$ with initial data $x$ and implementing the control $\alpha(\, \cdot\, )$. The functions $\ell: \R^3 \times A \times [0,+\infty) \to \R$ and $g: \R^3 \to \R$, called, respectively, the running cost and the final cost, are given bounded and continuous functions.

Given the dynamical system, it is crucial to define a suitable cost functional, in order to obtain the optimal strategy to control the epidemic. One can include several terms related to the sanitary impact and to other social effects. In particular, we take into account the following components:
\begin{enumerate}
\item social and economic costs due to the disease itself; this includes lost working hours, psychological consequences of isolation, etc.,
\item hospitalisation costs, including those related to Intensive Care Units (ICU),
\item  costs to implement the restrictive measures,
\item  vaccination costs.
\end{enumerate}
We assume that every cost component is an  increasing convex function, for both regularity and better representation of real costs (see \cite{OKS_2022}). For the socio-economic costs, we suppose that the infective are affected more than the others, therefore we take
$$ c_1 \ i^2 + \frac{c_1}{2} \ (1-i)^2,$$
where $c_1>0$ is an appropriate constant. For the second component, hospitalisations, we simply take
$$ c_2\ i^2, \quad c_2 > 0.$$ 
Since everyone is affected by social restrictions, we assume the third component to be proportional to the entire population, which is equal to $N=1$, so we simply get 
$$ c_\lambda \ \lambda^2, \quad c_\lambda >0. $$
Finally, for the vaccination costs, we suppose that there is a fixed cost, not proportional to the amount of people that are vaccinated, plus a second term that is indeed proportional to the susceptible, thus getting
$$ (c^0_\nu + c_\nu \ s^2)\ \nu^2, \quad c^0_\nu,\ c_\nu >0. $$
Summing all the components, we obtain the following running cost:
\begin{equation} \label{running_cost}
\ell(s,e,i,\lambda,\nu) = (c_1+c_2) \ i^2 + \frac{c_1}{2} \ (1-i)^2 +\ c_\lambda \ \lambda^2 + (c^0_\nu + c_\nu \ s^2)\ \nu^2.
\end{equation}
We point out that the choice of the constants $c_1,\ c_2,\ c_\lambda,\ c^0_\nu,\ c_\nu$ is a rather delicate task, as they should be fitted according to experimental data for the particular disease and population one is considering, but this goes beyond the scope of the present work. We also remark that what is most important in this control approach is not the precise value of each constant, but rather their relative magnitudes. For all the simulations in Section \ref{sec:numerical} we have used the values reported in following table:
\begin{center}
    \begin{tabular}{|c|c|}
    \hline
    \textbf{Constant} & \textbf{Value} \\
    \hline
    $c_1$ & 3.5 \\ 
    \hline
    $c_2$ & 14 \\ 
    \hline
    $c_\lambda$ & 0.35 \\ 
    \hline
    $c^0_\nu$ & 0.025 \\
    \hline
    $c_\nu$ & 0.05 \\ 
    \hline
    \end{tabular}
\end{center}
These values were chosen in order to prioritise the different cost components we considered. For the case of Covid-19, it has been estimated \cite{WFCJHCLY_2021,AP_2022,FSMTMPFB_2022} that hospitalisations are more expensive than lock-down policies, which are in turn way more expensive than vaccinations. In particular, the cost of a hospital bed in an intensive care unit is several orders of magnitude higher than a vaccine.

For the final cost, we can take into account the distance from a given desired state for each of the compartments. In our case, we want to penalise the scenarios in which the epidemic is not over by time $T$, therefore we set
\begin{equation} \label{final_cost}
g(s(T),e(T),i(T))= c_i \, (i(T)-\overline i)^2\ +\ c_e \, (e(T)-\overline e)^2,
\end{equation}
with $ c_i=c_e=10 \cdot c_1$ and $ \overline{i}=\overline{e}=0$.

We can also take into account a global bound for the Intensive Care Units (ICU), since their number is limited and it is one of the strong constraints for healthcare systems. In practice, this implies that the global number of people in ICU must stay always below that bound. One way to deal numerically with state constraints is described in Section \ref{sec:SLvsDAL}.

\section{The Dynamic Programming approach and semi-Lagrangian schemes}
\label{sec:dp}

In this section, we present the Dynamic Programming (DP) approach to optimal control problems and the related semi-Lagrangian scheme to obtain the numerical solution. We refer to \cite{FG_1999, BCD_2008, FF_SIAM2013, FF_2016} for a complete description of the method.
 
 We start by defining an auxiliary function, called the \textit{value function},
\begin{equation} \label{value}
v(x,t) := \inf_{\alpha \in \A} J_{x,t} (\alpha),
\end{equation}
which represents the best price we can pay starting from $x$ at time $t$. By means of the Dynamic Programming Principle \cite{B_1957, BCD_2008}, the value function, which is typically non-differentiable in optimal control problems, can be characterised as the unique viscosity solution {\cite{CL_1983,CEL_1984,E_1998} of the following Hamilton-Jacobi-Bellman equation for $x \in \R^3$ and $t \in (0,T)$:
\begin{equation} \label{HJB} 
\begin{cases}
- v_t (x,t) +\max\limits_{a \in A} \left\{ -f(x,a,t)\cdot \nabla v (x,t) -\ell(x,a,t) \right\}=0,\\
v(x,T)=g(x).
\end{cases}
\end{equation}
Viscosity solutions of partial differential equations, introduced by Crandall and Lions in the 1980s \cite{CL_1983,CEL_1984}, allow us to consider non-differentiable functions that satisfy \eqref{HJB} in a pointwise \textit{weak} sense. Nonetheless, unlike other notions of weak solutions, we can still get uniqueness results \cite{BCD_2008,FF_SIAM2013} for viscosity solutions of Hamilton-Jacobi equations, which are essential for the well-posedness of the problem.

Once the value function is known, optimal controls in feedback form can be computed by solving the following optimisation problem:
\begin{equation}\label{feedback}
a^* (x,t) = \argmax_{a\in A} \left\{ -f(x,a,t) \cdot \nabla v(x,t) -\ell(x,a,t) \right\}.
\end{equation}

Due to its non-linearity, equation \eqref{HJB} typically does not admit solutions in closed form, hence numerical schemes are needed to obtain an approximate solution. Constructing an approximation of a nonlinear partial differential equation has two main difficulties: the first is to deal with non regular solutions – typically just Lipschitz continuous – and the second is the high-dimensionality of the problem, since the number of equations can be rather high for compartmental models. We briefly describe the semi-Lagrangian method, which naturally follows the continuous control problem, via a discrete Dynamic Programming principle (see \cite{FF_SIAM2013,FF_2016} for details). First, we introduce a semi-discrete scheme with time step $\Delta t: = [(T-t)/n_{\max}]$, where $n_{\max}$ is the number of time steps, for $x \in \R^3$:
\begin{equation}\label{SL}
\begin{cases}
V^{n_{\max}}(x)=g(x), \\
V^{n}(x)=\min\limits_{a\in A}[V^{n+1}(\overline{x}) +\Delta t\,\ell(x, a, t_n)], \qquad n= n_{\max}-1,\ldots, 0,
\end{cases}
\end{equation}
where $V^n(x):=V(x, t_n)$, $t_n=t+n \, \Delta t$, $t_{n_{\max}} = T$ and $\overline{x}=x+\Delta t \, f(x, a, t_n)$. This is a backward problem, where we start from the final condition and we get back to the initial time. In order to obtain a fully discrete, explicit scheme, the term $V^{n+1}(\overline{x})$ is treated by interpolation on a grid, since $\overline{x}=x+\Delta t \, f(x, a, t_n)$ in general is not a grid point. This approximation leads to an a priori error estimate in terms of $\Delta t$ and $\Delta x$ that typically shows convergence of order $1/2$ to the value function, due to the fact that $v$ is just Lipschitz continuous (see \cite{FG_1999} for details). 

With the discrete value function at hand, still from \eqref{SL}, we can also synthesise the discrete feedback controls $a^*(x,t_n)$ and then reconstruct the optimal trajectories:
\begin{equation} \label{ricostruzione}
\begin{cases}
y^*(t_{n+1})=y^*(t_{n}) + \Delta t \, f(y^*(t_n),a^*(y^*(t_n),t_n),t_n), \\
y^*(t_0)=x.
\end{cases}
\end{equation}
Note in particular that $\alpha^*(t_n) := a^*(y^*(t_n),t_n)$ provides the corresponding discrete open-loop control.

It is worth noting that the semi-Lagrangian scheme is intrinsically parallel, since the computation of the value function on each grid node can be assigned to a single processor. Moreover, a single synchronisation among the processors is needed at each time step, but the computation of the new values only depends on the previous iteration and it requires the same amount of operations per node, including the location of the foot of the characteristic, the interpolation and the update. Finally, we point out that for our problem the structured grid and the positively invariant region make it possible to reduce the actual computation to a subset of the discrete state-space containing approximately $18\%$ of the nodes.


\section{The variational approach and the Direct-Adjoint Looping method}
\label{sec:pontryagin}

We now describe an alternative approach to Dynamic Programming, based on Pontryagin's maximum (or minimum) principle, for a finite horizon optimal control problem \cite{PBGM_1964,K_2004}. We keep the same notation adopted in the previous section. Instead of defining the value function, we want to find some conditions that make the cost functional $J_{x,t}(\alpha)$ stationary.

Assuming that the final cost $g$ is regular, the cost functional can be expressed as
\begin{equation}\label{costo_pontryagin}
J_{x,t}(\alpha) = \left[ \int_{t}^T \ell(y_{x,t} (\tau),\alpha(\tau),\tau) + \frac{d}{d\tau}\ g(y_{x,t}(\tau)) \ d\tau \right] + g(x).
\end{equation}
Since $t$ and $x=y_{x,t}(s) \big\lvert_{s=t}$ are fixed, minimising $J_{x,t}$ is equivalent to minimising
$$
\begin{aligned} 
\tilde{J}_{x,t}(\alpha) & = \int_{t}^T \ell(y_{x,t} (\tau),\alpha(\tau),\tau) + \frac{d}{d\tau}\ g(y_{x,t}(\tau)) \ d\tau = \\
&= \int_{t}^T \ell(y_{x,t} (\tau),\alpha(\tau),\tau) + \left[ \nabla g(y_{x,t}(\tau)) \right]\tpose \cdot \dot{y}_{x,t}(\tau)\ d\tau.
\end{aligned}
$$
Keeping in mind that along the optimal trajectories
$$  
f(y(\tau),\alpha(\tau),\tau) - \dot{y}(\tau) =0,
$$
we form the augmented cost functional
\begin{equation} \label{augmented}
J^a_{x,t}(\alpha) = \int_{t}^T \ell(y_{x,t} (\tau),\alpha(\tau),\tau) + \left[ \nabla g(y_{x,t}(\tau)) \right]\tpose \cdot \dot{y}_{x,t}(\tau) + p\tpose(\tau) \cdot \left[ f(y_{x,t}(\tau),\alpha(\tau),\tau) - \dot{y}_{x,t}(\tau) \right] \ d\tau
\end{equation}
by introducing the Lagrange multipliers $p(\tau)=\left( p_1(\tau),p_2(\tau),p_3(\tau) \right)$. From now on, to keep the notation lighter, we will omit the dependence of the solution from the initial data $\{x,t\}$ and all the explicit dependencies on time. Moreover, we assume for a moment that the controls are unconstrained. We then define the Hamiltonian
\begin{equation} \label{hamiltonian}
\mathcal{H}(y,\alpha,p,\tau) := \ell(y,\alpha,\tau) + p\tpose \cdot f(y,\alpha,\tau)
\end{equation}
and, by imposing stationarity on $J^a_{x,t}$ at $\alpha^*$, we get the following optimality system
\begin{equation} \label{pontryagin_var_zero}
\left\{
\begin{aligned}
&\dot{y^*}(\tau) = \frac{\partial \mathcal{H}}{\partial p} (y^*,\alpha^*,p^*,\tau) \\
&\dot{p^*}(\tau) = - \frac{\partial \mathcal{H}}{\partial y} (y^*,\alpha^*,p^*,\tau) \\
& \frac{\partial\mathcal{H}}{\partial \alpha} (y^*,\alpha^*,p^*,\tau) = 0,
\end{aligned}
\right.
\end{equation}
for all $\tau \in (t,T)$ and with mixed initial and final conditions
\begin{equation} \label{pontryagin_IFC}
\left\{
\begin{aligned}
& y^*(t)=x \\
& p^*(T)=\nabla_y \, g (y^*(T)).
\end{aligned}
\right.
\end{equation}
More explicitly, omitting all the dependencies for simplicity, we can write the following two-point boundary-value problem:
\begin{equation} \label{pontryagin_system}
\left\{
\begin{aligned}
&\dot{y^*} = f \\
&\dot{p^*} = - \nabla_y \, f \tpose \cdot p^* - \nabla_y \, \ell \\
&\nabla_\alpha \, f \tpose \cdot p^* + \nabla_\alpha \, \ell = 0 \\
& y^*(t)=x \\
& p^*(T)=\nabla_y \, g (y^*(T))
\end{aligned}
\right.
\end{equation}
which in our case is composed of eight equations: three for the state $y^*=(s^*,e^*,i^*)$, three for the co-state (or adjoint variables) $p^*$ and two for the optimality condition, plus the initial and final conditions.

We now briefly summarise the basic Direct-Adjoint Looping method to solve system \eqref{pontryagin_system} and we refer to \cite{K_2004} for further details: 
\begin{enumerate}
    \item discretise the time interval $[t,T]$ into $N$ subintervals, thus generating the discrete times $t_0=t, \ldots, t_N=T$, and set a descent step $\sigma < 1$ and a tolerance $\epsilon \ll 1$;
    \item select an initial discrete approximation of the control \\ $\alpha^{(0)}=(\alpha^{(0)}_0,\ldots,\alpha^{(0)}_N)\tpose$ and set $k=0$;
    \item once the control is fixed, integrate the forward equations for $y$ with initial data $y(t)=x$, obtaining a discrete approximation of the state $y^{(k)}=(y^{(k)}_0,\ldots,y^{(k)}_N)\tpose$;
    \item integrate backwards the equations for the co-state $p$ with final data 
    $$ p(T)=\nabla_y \, g (y^{(k)}_N), $$
    thus obtaining $p^{(k)}=(p^{(k)}_0,\ldots,p^{(k)}_N)\tpose$;
    \item if 
    $$  \bigg\lvert\bigg\lvert \, \frac{\partial \mathcal{H}^{(k)} }{\partial \alpha} \, \bigg\lvert\bigg\lvert := \bigg\lvert\bigg\lvert \, \nabla_\alpha \, f \tpose \cdot p^{(k)} + \nabla_\alpha \, \ell \, \bigg\lvert\bigg\lvert < \epsilon, $$ 
    terminate the procedure and the approximation of the optimal couple $(\alpha^*,y^*)$ is $(\alpha^{(k)},y^{(k)})$; otherwise, set 
    $$ \alpha^{(k+1)} = \alpha^{(k)} - \sigma \frac{\partial \mathcal{H}^{(k)} }{\partial \alpha}$$
    and repeat the procedure from step 3.
\end{enumerate}

Multiple techniques for nonlinear optimisation problems can be found in literature, in order to refine the algorithm and make it faster or more reliable, e.g. choosing the best descent step $\sigma^{(k)}$ at each iteration or using  inexact line search algorithms, see for instance \cite{NW_2006}. Moreover, to deal with constraints on the controls, the third equation in \eqref{pontryagin_var_zero} has to be replaced by the minimum principle
$$
\mathcal{H}(y^*(\tau),\alpha^*(\tau),p^*(\tau),\tau) \leq \mathcal{H}(y(\tau),\alpha(\tau),p(\tau),\tau)
\quad \mbox{for a.e. }\tau \in (t,T)\mbox{ and } \alpha \in \A\,,
$$
or rephrased as a variational inequality
$$
 \frac{\partial \mathcal{H}}{\partial \alpha}(y^*(\tau),\alpha^*(\tau),p^*(\tau),\tau) \cdot (a-\alpha^*(\tau))\ge 0 \qquad \mbox{for a.e. } \tau \in (t,T) \mbox{ and } a\in A, 
$$
and the optimality system can be solved using classical methods in nonlinear constrained optimisation \cite{multipliers-H,multipliers-P}. 

In our experiments, we always employ \textit{box constrained} controls, i.e. the control set $A$ is just the rectangle $[0,0.9] \times [0,\nu_{\max}]$.
Hence, we combine the update of the controls in step 5 of the algorithm with a projection step on the control set $A$, yielding 
$$ \alpha^{(k+1)} = \Pi_{A}\left\{ \alpha^{(k)} - \sigma \frac{\partial \mathcal{H}^{(k)} }{\partial \alpha}\right\}\,,$$
where, recalling that $\alpha=(\lambda,\nu)$, the projection $\Pi_A=(\Pi_\lambda,\Pi_\nu)$ acts component-wise, respectively as $\Pi_\lambda\{\,\cdot\,\}=\max\{0,\,\min\{\,\cdot\,,0.9\}\}$ and 
$\Pi_\nu\{\,\cdot\,\}=\max\{0, \min \{ \,\cdot\, ,\,\nu_{\max} \}\}$. In addition, we use a simple bisection method on the descent step $\sigma$ in order to guarantee a strictly monotone decrease in the cost functional,    
while checking the convergence of the optimisation procedure using the following criterion: 
$$
\big\lvert J_{x,t}(\alpha^{(k+1)}) - J_{x,t}(\alpha^{(k)}) \big\lvert < \epsilon.
$$

We finally observe that, in the present setting, the above variational inequality for $\alpha^*(\cdot)=(\lambda^*(\cdot),\nu^*(\cdot))$ simplifies, component-wise and for a.e. $\tau \in (t,T)$, into
\begin{subequations} \label{firstOrderConditions}
\begin{equation}
\frac{\partial \mathcal{H}}{\partial \lambda}(y^*(\tau),\alpha^*(\tau),p^*(\tau),\tau)
\begin{cases}
 \, =0 & \mbox{if } 0<\lambda^*(\tau)<0.9 \\
 \, \ge 0 & \mbox{if } \lambda^*(\tau)=0 \\
 \, \le 0 & \mbox{if } \lambda^*(\tau)=0.9  
\end{cases}
\end{equation}
and
\begin{equation}
\frac{\partial \mathcal{H}}{\partial \nu}(y^*(\tau),\alpha^*(\tau),p^*(\tau),\tau)
\begin{cases}
 \, =0 & \mbox{if } 0<\nu^*(\tau)<\nu_{\max} \\
 \, \ge 0 & \mbox{if } \nu^*(\tau)=0 \\
 \, \le 0 & \mbox{if } \nu^*(\tau)=\nu_{\max}
\end{cases}
.
\end{equation}
\end{subequations} 
This, together with the positive semi-definiteness of the Hessian of $\mathcal{H}$ (Legendre condition), provides first and second order necessary conditions for an optimal control to be a local minimum, and it can be readily checked just after the convergence of the algorithm. More precisely, we verify that, for a.e. $\tau \in (t,T),$
\begin{subequations} \label{secondOrderCondition}
\begin{alignat}{3}
\frac{\partial^2 \mathcal{H}}{\partial \alpha^2}(y^*(\tau),\alpha^*(\tau),p^*(\tau),\tau)\geq 0 \quad &\text{ if } \alpha^*(\tau) \text{ is not on the boundary of } A, \label{secondOrderCondition_a} \\ 
\frac{\partial^2\mathcal{H}}{\partial\lambda^2} (y^*(\tau),\alpha^*(\tau),p^*(\tau),\tau)\geq 0 \quad &\text{ if } \lambda^*(\tau) \text{ is internal and } \nu^*(\tau) \text{ on the boundary,} \label{secondOrderCondition_b} \\
\frac{\partial^2\mathcal{H}}{\partial\nu^2} (y^*(\tau),\alpha^*(\tau),p^*(\tau),\tau)\geq 0 \quad &\text{ if } \nu^*(\tau) \text{ is internal and } \lambda^*(\tau) \text{ on the boundary.} \label{secondOrderCondition_c}
\end{alignat}
\end{subequations}


\section{Comparing and combining semi-Lagrangian schemes and the Direct-Adjoint Looping method}
\label{sec:SLvsDAL}

In this section, we briefly discuss the features of the two approaches we presented above, in particular their pros and cons, which depend on the problem one wants to solve. Moreover, we propose a simple but effective way to couple them in order to produce reliable solutions. The key point is to exploit the theoretical convergence of the semi-Lagrangian scheme to the global minimum of the cost functional.

\subsection{Feedback vs. open loop}
The major perk of the Dynamic Programming approach is that it provides feedback controls, i.e. controls that are a function of the state of the system. This means that the control policy can instantly react to small perturbations or uncertainty in the data. The variational approach, instead, leads to open loop controls, i.e. controls that are only a function of time. As a consequence, they can be no longer optimal if there are, for instance, some model errors or external disturbances. 

\subsection{The role of initial data}
In the Dynamic Programming approach, the variables of the value function are the initial data $(x,t)$ of the dynamical system. This implies that if the initial data change, we do not have to compute the value function again, but only the optimal trajectories. This allows us to build a static controller that takes the initial data as an input and outputs the control policy in \eqref{ricostruzione}. This property can be useful in case different strategies have to be computed, starting from various initial data. The DAL algorithm, instead, must be run again completely if we change the initial data of the dynamical system.

\subsection{Convergence and error estimates}
For semi-Lagrangian schemes, theoretical results show convergence to the optimal solution, i.e. to the \textit{global} minimum of the cost functional, with an a priori error estimate of order $1/2$ in the case of a Lipschitz continuous  value function, depending on the time and space discretisations \cite{FG_1999,FF_SIAM2013}. On the contrary, the variational approach only uses necessary conditions for a control to be optimal, therefore we have no guarantee that the iterative procedure will converge to the global minimum. Unless further hypotheses on $J_{x,t}$ are verified, such as convexity or coercivity, it may converge to different stationary points depending on the initial choice for $\alpha^{(0)}$, $\sigma$ and $\epsilon$. It may even not converge at all for some choices of the initial parameters. Nevertheless, when it converges, the accuracy is much higher than that of a semi-Lagrangian scheme, since only the temporal interval is discretised; both the state-space and the control set are treated as continuous. Moreover, there is no interpolation error in the DAL algorithm: the only source of error is the integration of the ordinary differential equations, which can be treated with a high-order numerical scheme. In the simulations presented in Section \ref{sec:numerical}, in fact, we will observe some differences in the control policies given by the two methods. In particular, those obtained with DAL appear smoother and their cost, computed a posteriori, is slightly smaller. Nevertheless, the difference in the trajectories is always compatible with the magnitude of the discretisation steps.

\subsection{State constraints}
\label{subsec:constraints}
The need to impose some state constraints in an optimal control problem is not uncommon: one may need the optimal trajectories to stay inside a certain subset $\Omega$ of the state-space. This implies substituting the space of admissible controls $\A$ with its subset
$$
\A_\Omega=\left\{ \alpha \in \A\ |\ y_{x,t} (\, \cdot \, ; \alpha): [t,T] \to \Omega \right\}.
$$
Numerically, we can force the trajectories to stay inside $\Omega$ (or, conversely, drop all the control strategies that let them outside $\Omega$) by simply changing the running cost $\ell(y,\alpha,\tau)$ so that it rapidly increases as soon as the trajectories $y$ leave $\Omega$.

\subsection{Computational cost}
The main difference between the two methods lays in their computational cost. Although the total number of iterations of the DAL algorithm cannot be known a priori, each iteration is going to be rather fast, since it mainly consists in integrating a system of ordinary differential equations. On the other hand, we can set the total number of temporal steps for the semi-Lagrangian scheme, but each iteration is going to be very expensive, since we need to build a local interpolation operator for every single node in the discrete state-space, for every single discrete control. As a result, the computational cost of semi-Lagrangian schemes drastically grows with the dimension of the state-space, i.e. with the number of equations in the dynamical system. This phenomenon is known as the \textit{curse of dimensionality}, expression coined by Bellman himself \cite{B_1957}, and makes serial implementations of SL schemes unfeasible for systems with more than 4 or 5 equations. However, in order to mitigate the computational efforts, several acceleration methods for particular Hamilton-Jacobi equations have been developed, such as Fast Marching \cite{S_1999} and Fast Sweeping methods \cite{FLZ_2009}, domain decomposition methods (for example \cite{CCFP_2012}) or schemes based on unstructured grids \cite{AFS_2020,CF_2021,AFS_2022,FKS_2023}.

\subsection{Combination of the two approaches}
The theoretical connection between Dynamic Programming and Pontryagin's principle has been widely investigated in literature. In particular, it is well known that, under smoothness assumptions, the gradient of the value function coincides with the adjoint variable along the optimal trajectory. In other words, the solutions of the optimality system \eqref{pontryagin_system} correspond to the characteristic curves of the HJB equation \eqref{HJB}, thus providing a representation formula for its solution. This result has been generalized in different directions to handle the non-smooth case, e.g. employing the notion of viscosity solutions or minmax solutions (see \cite{BCD_2008,AKK_2021,YD_2021} and references therein). From a numerical point of view, this connection has been exploited in different contexts and with different goals. In \cite{CM_2010}, under the assumption that the optimal control in \eqref{pontryagin_var_zero} can be expressed as a function of the state and the adjoint variables only, the optimality system \eqref{pontryagin_system} can be rewritten as a Cauchy problem for $(y(\cdot),p(\cdot))$ with a forward-forward structure, where the initial value $y(t)=x$ is given, while $p(t)$ is an unknown. The system is then solved by a Newton-based shooting method to match the final condition $p(T)=\nabla_y\, g(y(T))$, using the gradient of the value function at $(x,t)$ as a good initial guess for $p(0)$. On the other hand, under similar assumptions in the smooth case with nonlinear, control-affine dynamics and unconstrained controls, in \cite{AKK_2021} the optimality system \eqref{pontryagin_system} is solved keeping the forward-backward structure and combined with a gradient descent method to minimize the cost functional of the problem. This provides a fast and accurate approximation of the value function and its gradient at single points in the state space without solving the HJB equation. These data are then employed as a training set for a polynomial model for the value function on the whole state space, thus remarkably reducing the complexity in the case of high-dimensional control problems.

We now propose a combination of the semi-Lagrangian scheme and the DAL method, which shares similarities with both techniques discussed above. 
We first use SL to compute offline an approximation of the value function on a reasonably fine grid on the whole state space. Then, given the initial data $y(t)=x$, we repeatedly synthesise online an optimal feedback control, using \eqref{feedback} while integrating the forward dynamics in \eqref{pontryagin_system} up to the final time. This provides the initial guess for the control in the DAL method. Note that here the control is defined for all the discrete times, not only at the initial time $t$, as for the adjoint variable $p(t)$ in \cite{CM_2010}. From this point on, the DAL method proceeds until convergence, using the forward-backward structure of the optimality system to reduce the cost functional as in \cite{AKK_2021}. Nevertheless, we observe that we do not need any further assumptions on the optimal control and on the cost functional to produce a descent direction: it is obtained by using the current control and the updated state/adjoint variables to evaluate the gradient of the Hamiltonian with respect to the control, possibly using a projection to easily handle the case of box constrained controls. Finally, we remark as well that the choice of nonlinear, control-affine dynamics is more related to the SEIR models considered in this paper than to the proposed resolution method. As long as the gradient of the dynamics with respect to the control can be computed analytically, the DAL method can work out of the box without any modification.


\section{Numerical simulations}
\label{sec:numerical}

In this section, we report the numerical results we obtained applying the two approaches and their combination to some variations of model \eqref{seir_controllato_ridotto}. The Direct-Adjoint Looping algorithm was implemented in C++ and run on a $1.4$ GHz Intel Core i5 quad-core CPU, whereas a parallel version of semi-Lagrangian scheme was implemented in CUDA C++ and run on a single NVIDIA GeForce RTX 2070 GPU. The temporal interval we consider is 3 years and the time unit is one trimester, therefore, in the notation adopted above, $t=0$ and $T=12$. The time step is $\Delta t = 0.05$, meaning 600 temporal iterations for both methods. For the the semi-Lagrangian scheme, we take a uniform 3D mesh in the state-space $[0,1]^3$ made of $150^3$ nodes.
To set some reasonable parameters, we consider an ideal infectious disease for which the latency period lasts on average $1/\varepsilon=10$ days and the mean infectious period is $1/\gamma=3$ weeks, hence $\varepsilon=9$, $\gamma=4$. In addition, we set
\begin{equation}
\label{beta}
\beta(\tau)=
\begin{cases}
4 \quad &\text{if } 2 \leq \overline{t} \leq 3, \\
16 \quad &\text{otherwise,}
\end{cases}
\end{equation}
where $\overline{t} \equiv \tau$ mod $4$, so that $R_t=R_0=4$  all year long except for one trimester, where it drops to 1. This choice of $\beta$ is to represent the natural decrease in the transmission rate that we observe during summer time for many infectious diseases, like chickenpox, influenza and even Covid-19 \cite{LY_1973,LS_2014,LHLZWHY_2021}. Moreover, we consider a vaccine with $p=90\%$ efficacy, meaning we substitute $\nu$ with $p\cdot\nu$ in \eqref{seir_controllato_ridotto} wherever it appears. Finally, we assume that there are 1000 infective, 3000 exposed and no recovered at time $t=0$. We normalise these quantities with the total Italian population $N=58,983,122$ \cite{ISTAT}, obtaining the initial data 
\begin{align*}
    &e_0=\frac{3000}{N} = 5.09 \cdot 10^{-5}, \\
    &i_0= \frac{1000}{N} = 1.70 \cdot 10^{-5}, \\
    &r_0 = 0, \\
    &s_0=1-e_0-i_0.
\end{align*}
The value for $s_0$ is calculated keeping in mind that $s_0+e_0+i_0+r_0=1$, since the four variables represent the fractions of population in each compartment.

Regarding the controls, we set 
$$ A = A(\tau) :=[0,0.9]\times [0,\nu_{\max}(\tau)],$$ 
where
\begin{equation} \label{numax}
\nu_{\max}(\tau)=
\begin{cases}
0 \quad &\text{if } \tau < 4, \\
(\tau-4) \quad &\text{if } 4 \leq \tau <5, \\
1 \quad &\text{if } \tau \geq 5,
\end{cases}
\end{equation}
in order to mimic the initial unavailability of vaccines when a new virus emerges. Although it is theoretically possible, we do not allow $\lambda$ to reach $1$, as it would mean preventing any contact between susceptible and infectious individuals and that is impossible to achieve, even with the strongest lock-down policies adopted for the Covid-19 pandemic. For the SL scheme, the control set $A$ is discretised with a uniform 2D mesh composed of $150^2$ discrete controls.

From now on, the superscript \textit{SL} will indicate the optimal solutions obtained with the semi-Lagrangian scheme and the reconstruction \eqref{ricostruzione}, whereas \textit{DAL} will be assigned to the results of the Direct-Adjoint Looping algorithm. The superscript \textit{SL-DAL} indicates the results of the combined scheme. For comparison purposes, in each test we evaluate the corresponding cost functional $J_{x,t}$ in \eqref{costo_orizz_fin} along the optimal trajectories computed by the two algorithms, discretised with a simple rectangular quadrature rule using the same step $\Delta t$ for the time interval $[0,T]$.

\subsection{Test 1: basic model}

We begin with the simplest case, that is the basic model \eqref{seir_controllato_ridotto} with running cost \eqref{running_cost} and final cost \eqref{final_cost}, meaning that our final goal is to end the infection by time $\tau=T$. We repeat the DAL algorithm with various initial guesses for $\lambda^{(0)}$ and $\nu^{(0)}$, including $\lambda^{SL}$ and $\nu^{SL}$, always obtaining the same results, which are reported in Figure \ref{fig:sim1}.

\begin{figure}[h]
    \centering
    \includegraphics[width=7cm]{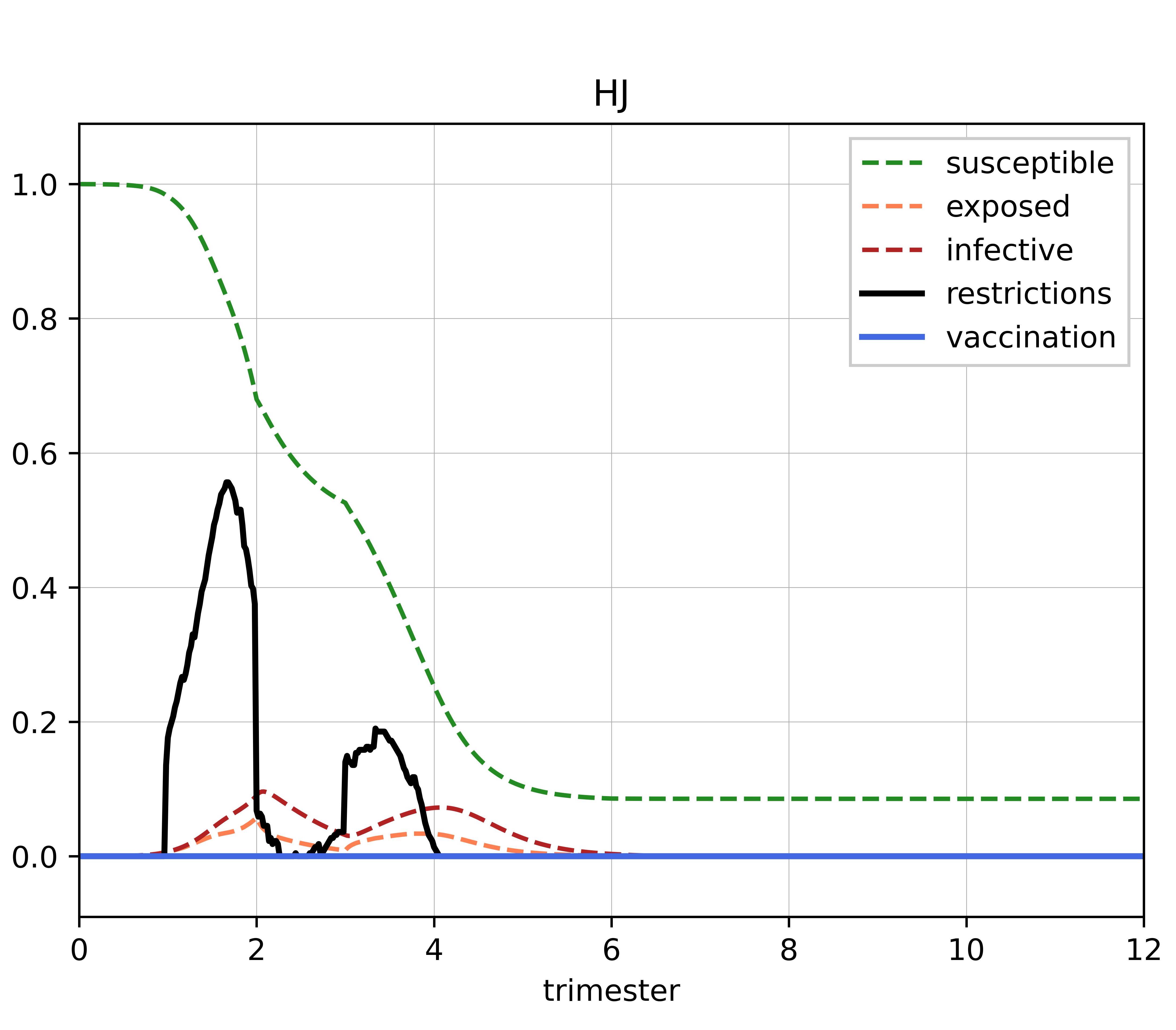}
    \includegraphics[width=7cm]{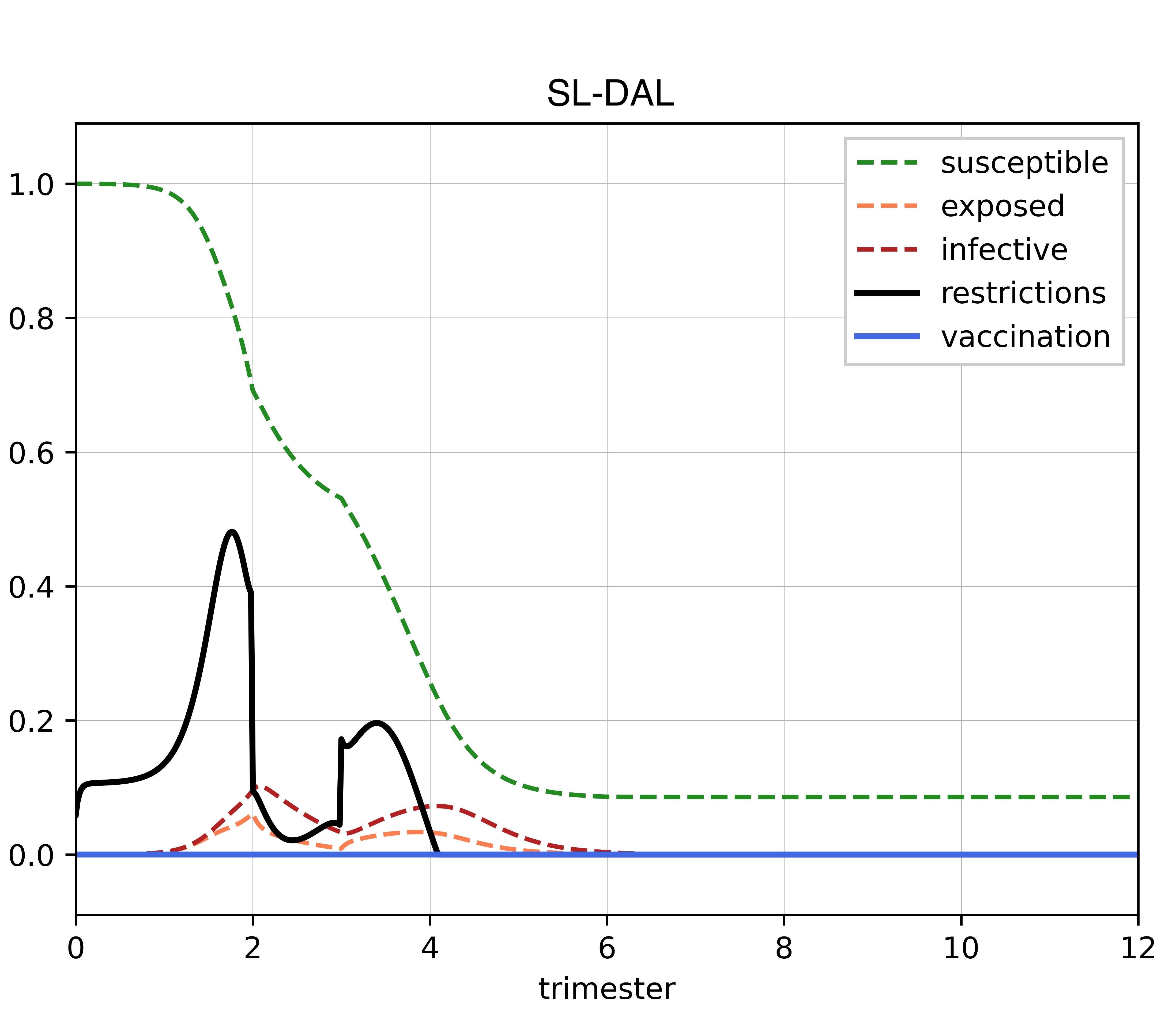}
    \caption{\textit{Test 1.} Optimal trajectories (dashed lines) and controls (full lines) for system \eqref{seir_controllato_ridotto}. Solutions obtained with the semi-Lagrangian scheme on the left, with Direct-Adjoint Looping on the right.}
    \label{fig:sim1}
\end{figure}

As already anticipated, the controls we get from the two methods are slightly different, but qualitatively the same. In particular, the optimal trajectories (dashed lines in Figure \ref{fig:sim1}) are nearly indistinguishable. We observe that the typical peak of infective (see Figure \ref{fig:grafici_seir} for reference) is cut using mobility restrictions, then, when a second peak is about to form, some milder restrictions are applied in order to bring the infective down to zero. Vaccination is never applied with this model.

\subsection{Test 2: temporary immunity}

For the second simulation we consider a scenario in which immunity – whether acquired by infection or vaccination – is not permanent. Suppose its mean duration is $1/\mu=9$ months. This means we have to add a path from $R$ back to $S$ in the diagram in Figure \ref{fig:diagramma_seir}, leading to the following dynamical system:
\begin{equation} 
\begin{cases}
s' = - \beta (1-\lambda) s i - \nu s + \mu r \\
e' = \beta  (1-\lambda) s i - \varepsilon e \\
i' = \varepsilon e - \gamma i \\
r' =  \gamma i + \nu s - \mu r
\end{cases}
\end{equation}
with the same parameters and initial data as before, plus $\mu=1/3$. Now the fourth equation is not independent anymore, since $r$ appears in the first equation, but we can still use the conservation property to express $r(\tau)=1-[s(\tau)+e(\tau)+i(\tau)]$ and thus reduce the system to
\begin{equation} \label{seir_immunità_temporanea}
\begin{cases}
s' = - \beta (1-\lambda) s i - \nu s + \mu\, [1- s - e - i] \\
e' = \beta  (1-\lambda) s i - \varepsilon e \\
i' = \varepsilon e - \gamma i
\end{cases}
\end{equation}
As before, we run the DAL algorithm with various initial guesses for $\lambda^{(0)}$ and $\nu^{(0)}$, including $\lambda^{SL}$ and $\nu^{SL}$. In this case as well, the results are always the same and are shown in Figure \ref{fig:sim2}.
\begin{figure}[h]
    \centering
    \includegraphics[width=7cm]{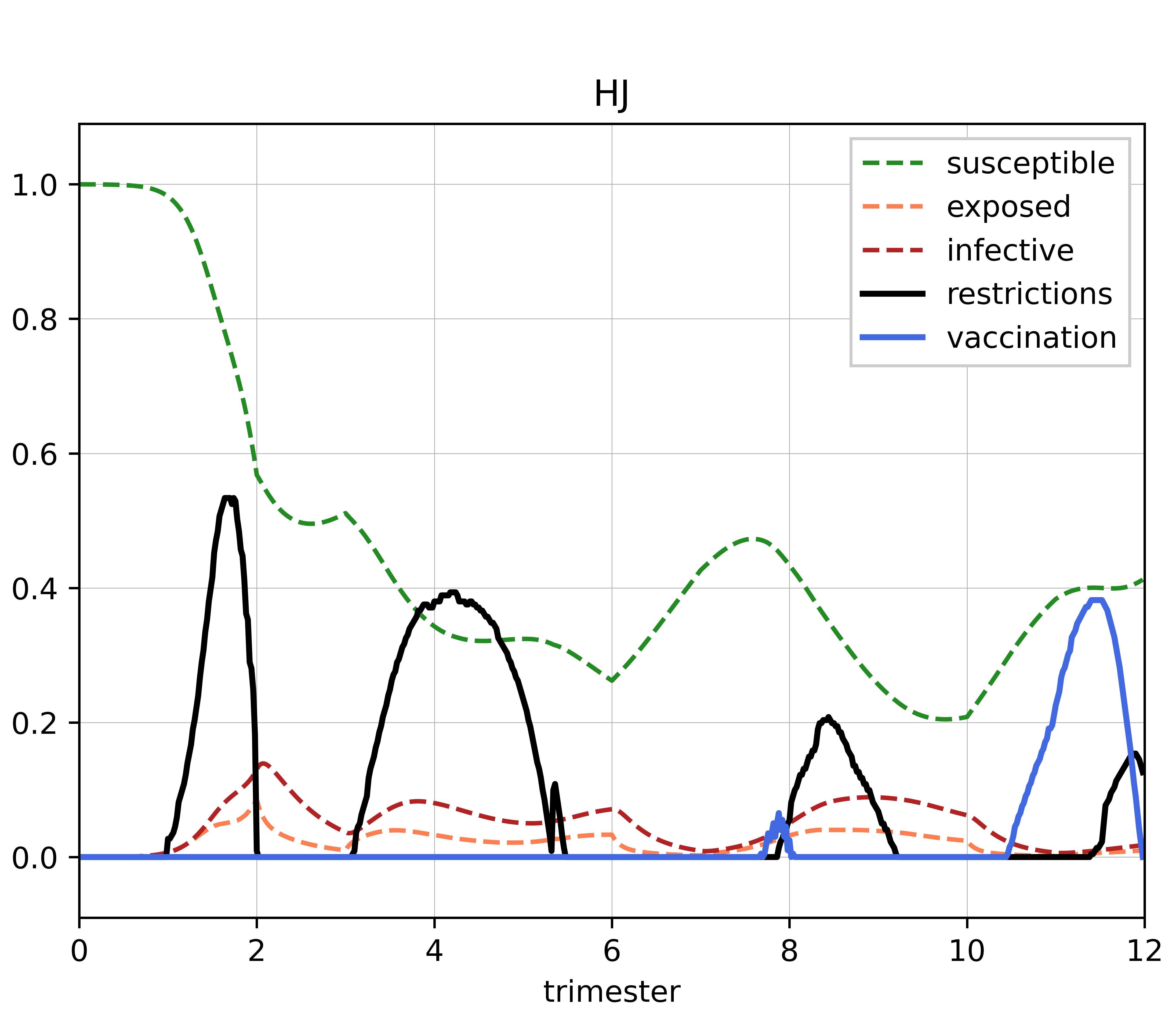}
    \includegraphics[width=7cm]{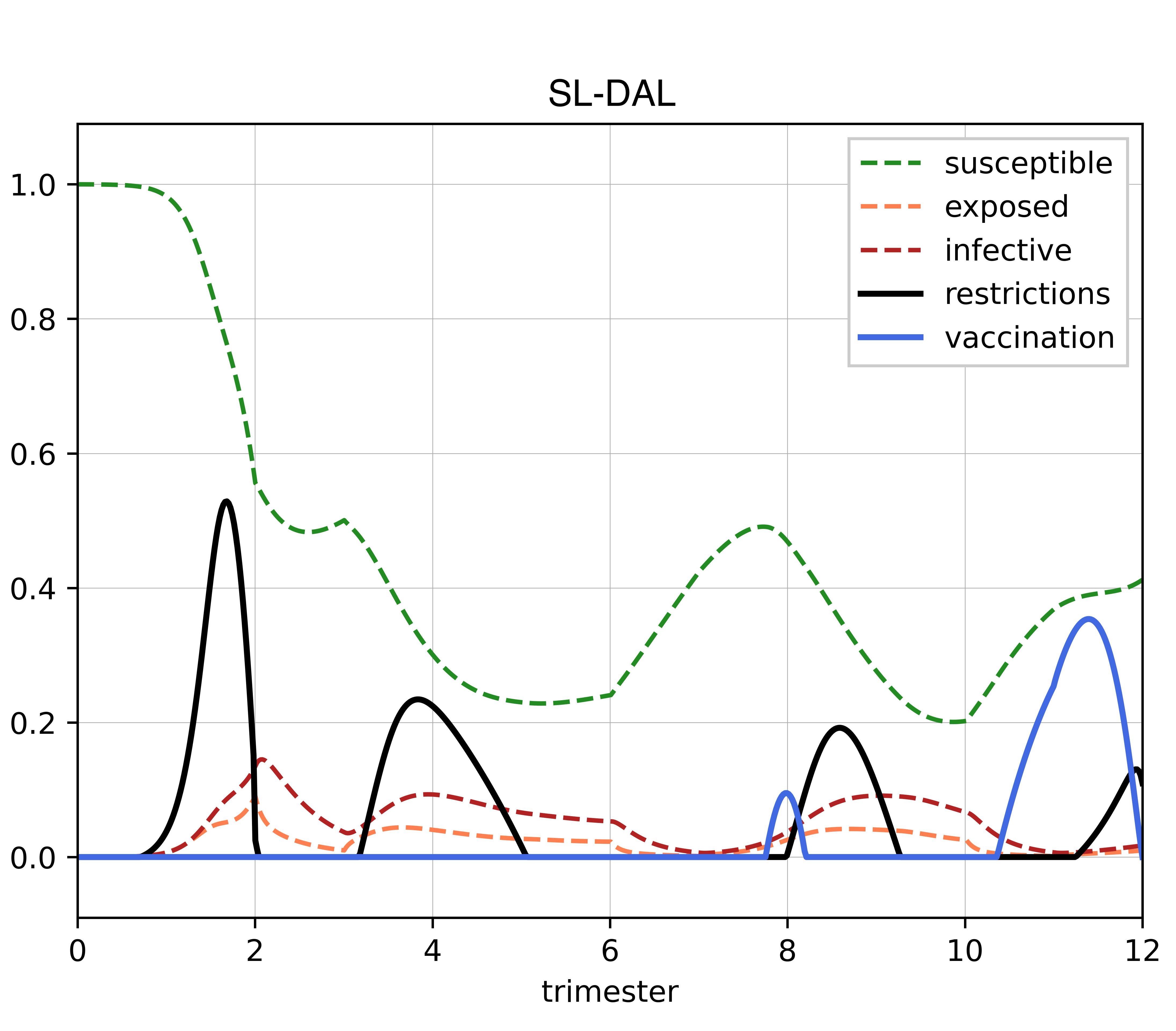}
    \caption{\textit{Test 2.} Optimal trajectories (dashed lines) and controls (full lines) for system \eqref{seir_immunità_temporanea}. Solutions obtained with the semi-Lagrangian scheme on the left, with Direct-Adjoint Looping on the right.}
    \label{fig:sim2}
\end{figure}
The two numerical schemes give slightly different results, but qualitatively the same. What is peculiar of models with temporary immunity is the oscillations in the trajectories due to the different rates with which the susceptible and recovered compartments exchange individuals, which are also reflected in the control policy. As a matter of fact, we can observe four different periods of mobility restrictions, the first between the first and second trimester and the last at the very end of the 3-year window, together with two waves of vaccination, a first small one around the eighth trimester and a bigger one at the end. The restrictions are milder and milder with the passing of time, while the vaccination rate increases.

\subsection{Test 3: border control}
For this simulation, we go back to system \eqref{seir_controllato_ridotto}, but we introduce a constant source term $\delta$ in the model, representing incoming people from the outside. Let $\delta_j$, $j=1,\ldots,4$, be the fractions of $\delta$ going respectively into the $s$, $e$, $i$, $r$ compartments. In order to show that restrictive measures and vaccines are not the only possible controls that one can apply to epidemiological models, for this simulation we will not consider vaccinations; we keep the restrictive measures $\lambda$ as before and introduce a new control $\phi$ representing the opening of borders, so that we can control the influx of individuals. We want $\phi=0$ to correspond to closing all the borders and $\phi=1$ to be equivalent to the uncontrolled scenario. A diagram of this variation of the SEIR model is represented in Figure \ref{fig:diagramma_seir_frontiere}.

\begin{figure}[!h]
\centering
\includegraphics[width=9cm]{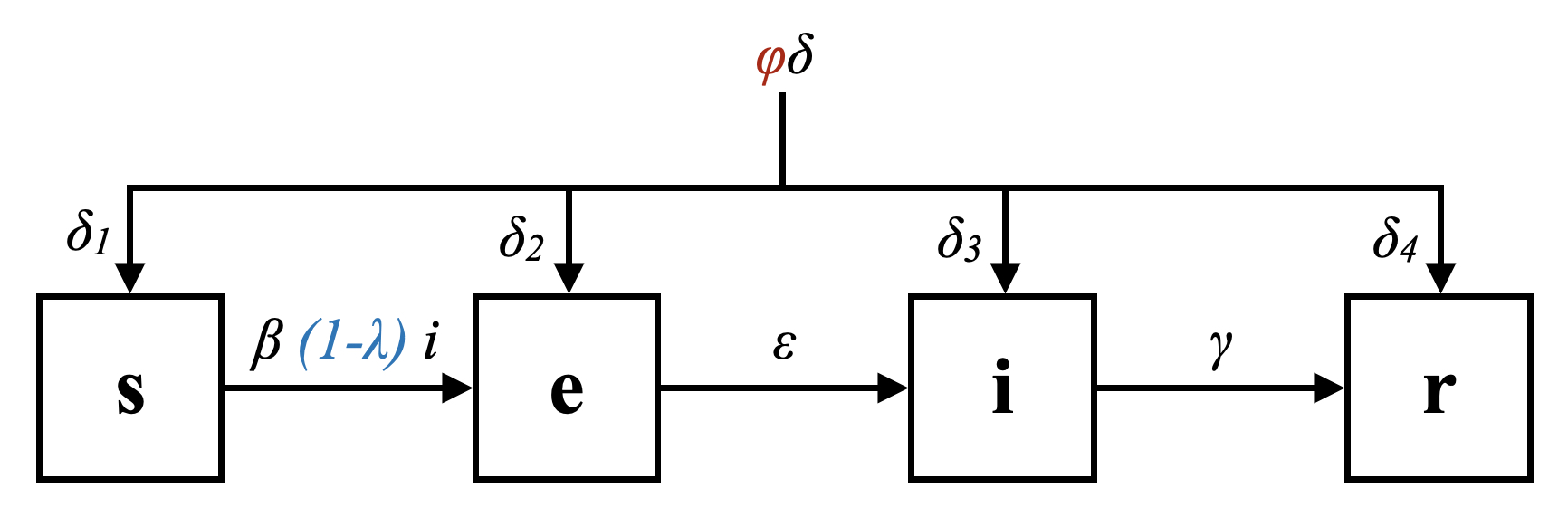}
\caption{Diagram of the normalised epidemic SEIR model with restrictions $\lambda$ and a control $\phi$ over incoming individuals.}
\label{fig:diagramma_seir_frontiere}
\end{figure}
In this case the controlled differential system becomes
\begin{equation} \label{seir_frontiere}
\begin{cases}
s' = - \beta (1-\lambda) s i + \phi \delta_1\\
e' = \beta  (1-\lambda) s i - \varepsilon e + \phi \delta_2\\
i' = \varepsilon e - \gamma i + \phi \delta_3\\
r' =  \gamma i + \phi \delta_4
\end{cases}
\end{equation}
Note that the conservation of population does not hold anymore, since 
\begin{equation}
    \label{no_conservazione}
    \frac{d}{d\tau} \left[s(\tau)+e(\tau)+i(\tau)+r(\tau)\right] = \delta \phi(\tau) \geq 0.
\end{equation}
However, the variable $r$ still does not appear in the first three equations, so it is always possible to solve the reduced system first and then eventually get $r(\tau)$ by integration of the fourth equation.

For the optimal control problem, we have to modify the running cost in order to:
\begin{itemize}
    \item eliminate the vaccination cost, since vaccines are not considered,
    \item take into account that the population is not constant,
    \item penalise closing the borders.
\end{itemize}
Therefore, we take out the term $(c^0_\nu + c_\nu \ s^2)\ \nu^2$ from the running cost, modify the restrictions cost to $c_\lambda \lambda^2 N(\tau)$ add a new term $c_\phi (1-\phi)^2 N(\tau)$, where $N(\tau)$ is the total population at time $\tau$. The reason for this is that we want the cost components relative to $\lambda$ and $\phi$ to be proportional to the entire population, since everyone is affected by the closure of the borders. The issue is that now, as already mentioned, the sum of the four variables is not equal to 1 for all $\tau \geq 0$, but only at $\tau=t=0$. From \eqref{no_conservazione} and the initial data we get 
\begin{equation}
    \label{popolazione}
    N(\tau)=1+ \delta \int_0^\tau \phi(\xi)\ d\xi,
\end{equation}
which we approximate as $ N(\tau) \approx 1 + \delta\, \tau \, \phi(\tau)$, obtaining the following running cost:
\begin{equation} \label{running_cost_frontiere}
 \ell(s,e,i,\lambda,\phi,\tau) = (c_1 + c_2) \, i^2 + \frac{c_1}{2} \, (1-i)^2 + c_\lambda \, \lambda^2 (1 + \delta \tau \phi) + c_\phi \, (1-\phi)^2 (1 + \delta \tau \phi).
\end{equation} 
Note that now $\ell$ explicitly depends on time. Moreover, we set the final cost \\ $g(s(T),e(T),i(T))=0$.

The constants relative to the incoming individuals that we used in this simulation are reported in the table below, while all the others are the same as the previous tests.
\begin{center}
    \begin{tabular}{|c|c|}
    \hline
    \textbf{Constant} & \textbf{Value} \\
    \hline
    $\delta$ & 0.75 \\
    \hline
    $\delta_1$ & $0.5 \cdot \delta$\\
    \hline
    $\delta_2$ & $0.01 \cdot \delta$  \\
    \hline
    $\delta_3$ & $0.005 \cdot \delta$ \\
    \hline
    $\delta_4$ & $0.485 \cdot \delta$ \\
    \hline
    $c_\phi$ & $0.15$ \\
    \hline
    \end{tabular}
\end{center}

\begin{figure}[h]
    \centering
    \includegraphics[width=7cm]{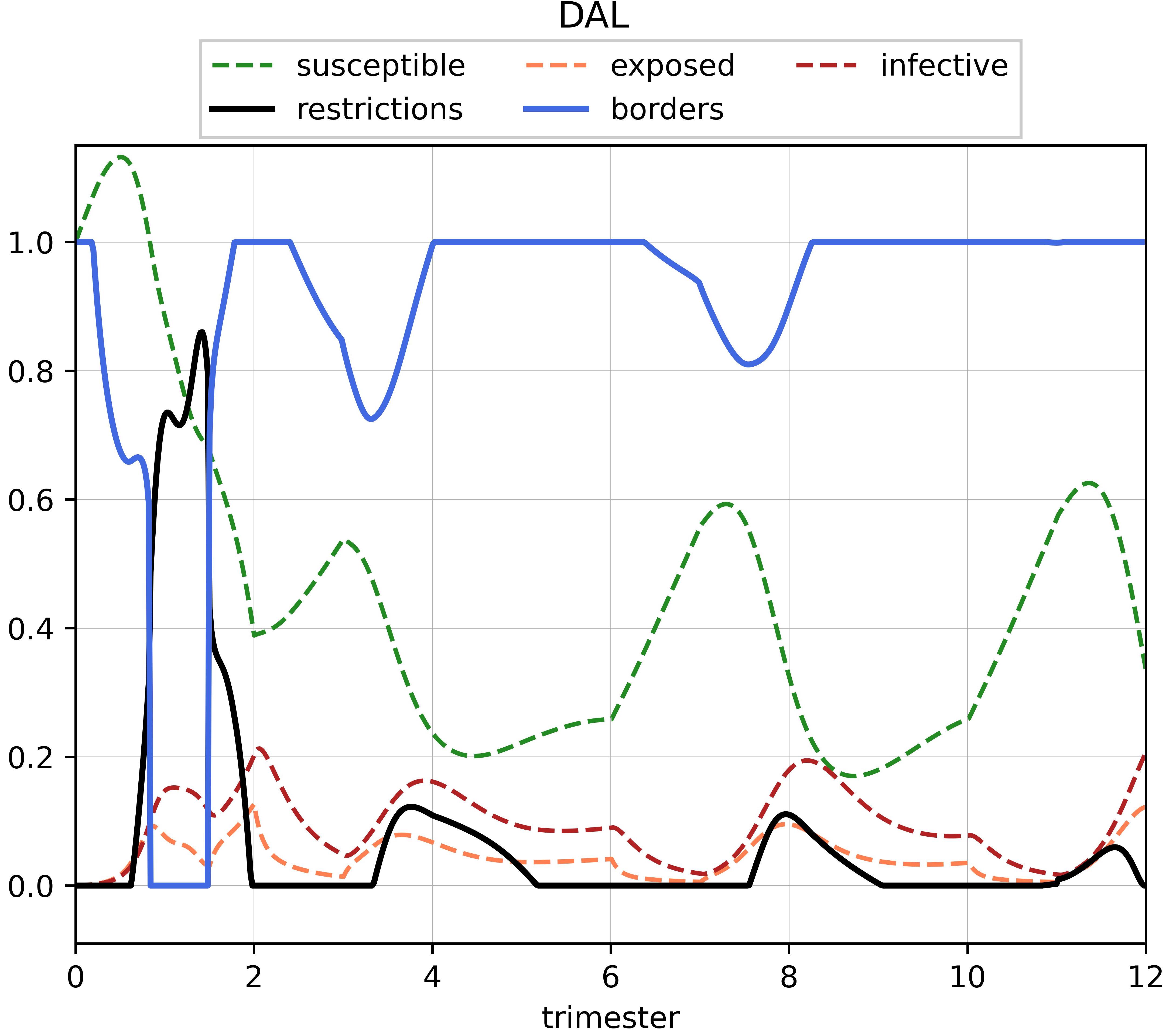}
    \includegraphics[width=7cm]{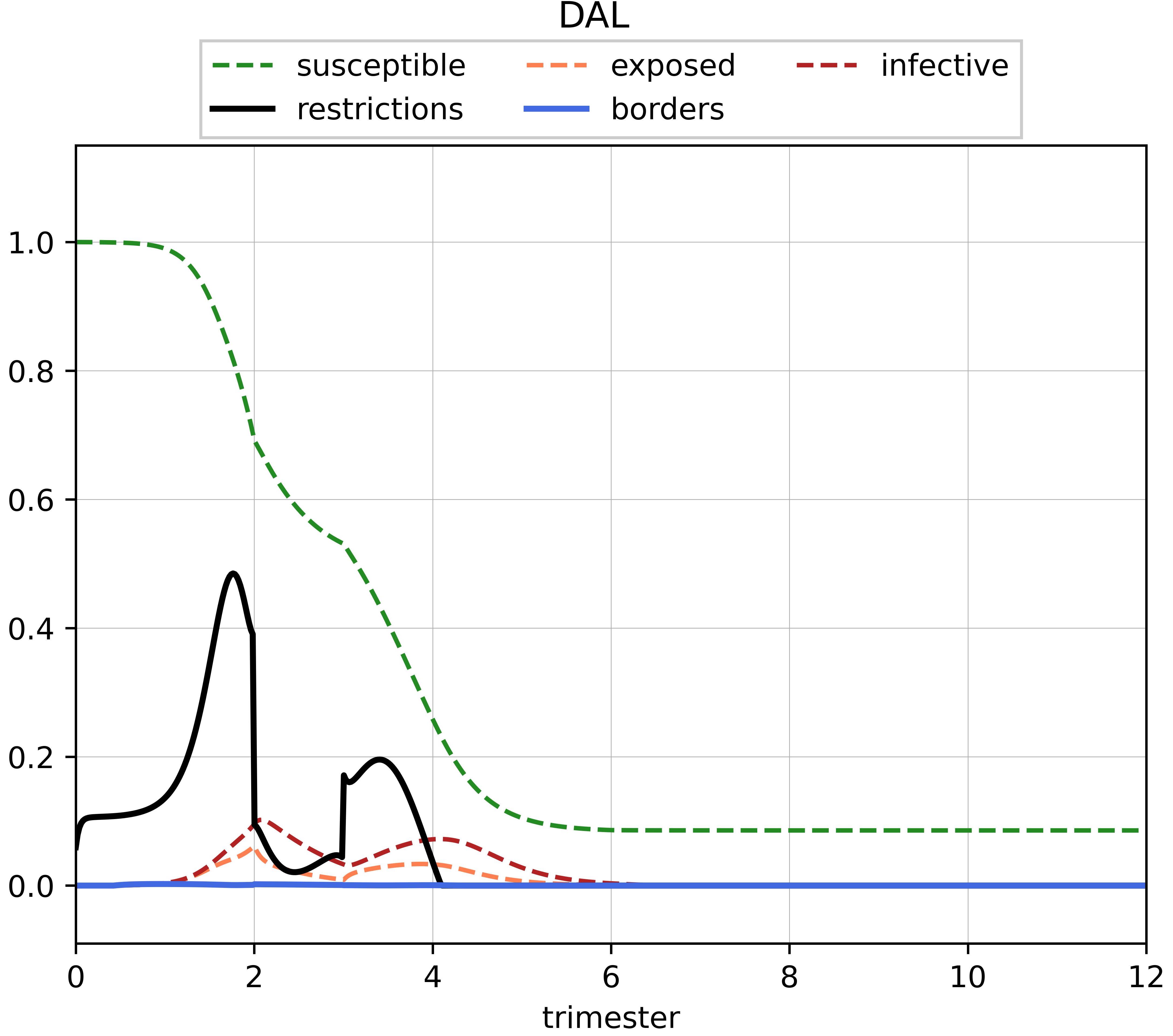}
    \caption{\textit{Test 3.} Optimal trajectories (dashed lines) and controls (full lines) for system \eqref{seir_frontiere} obtained with the DAL method. On the left with initial guess $\lambda^{(0)}=[0,\ldots,0]$ and $\phi^{(0)}=[1,\ldots,1]$, and on the right with initial guess $\lambda^{(0)}=\phi^{(0)}=[0,\ldots,0]$.}
    \label{fig:sim3wrong}
\end{figure}

\begin{figure}[h]
    \centering
    \includegraphics[width=7cm]{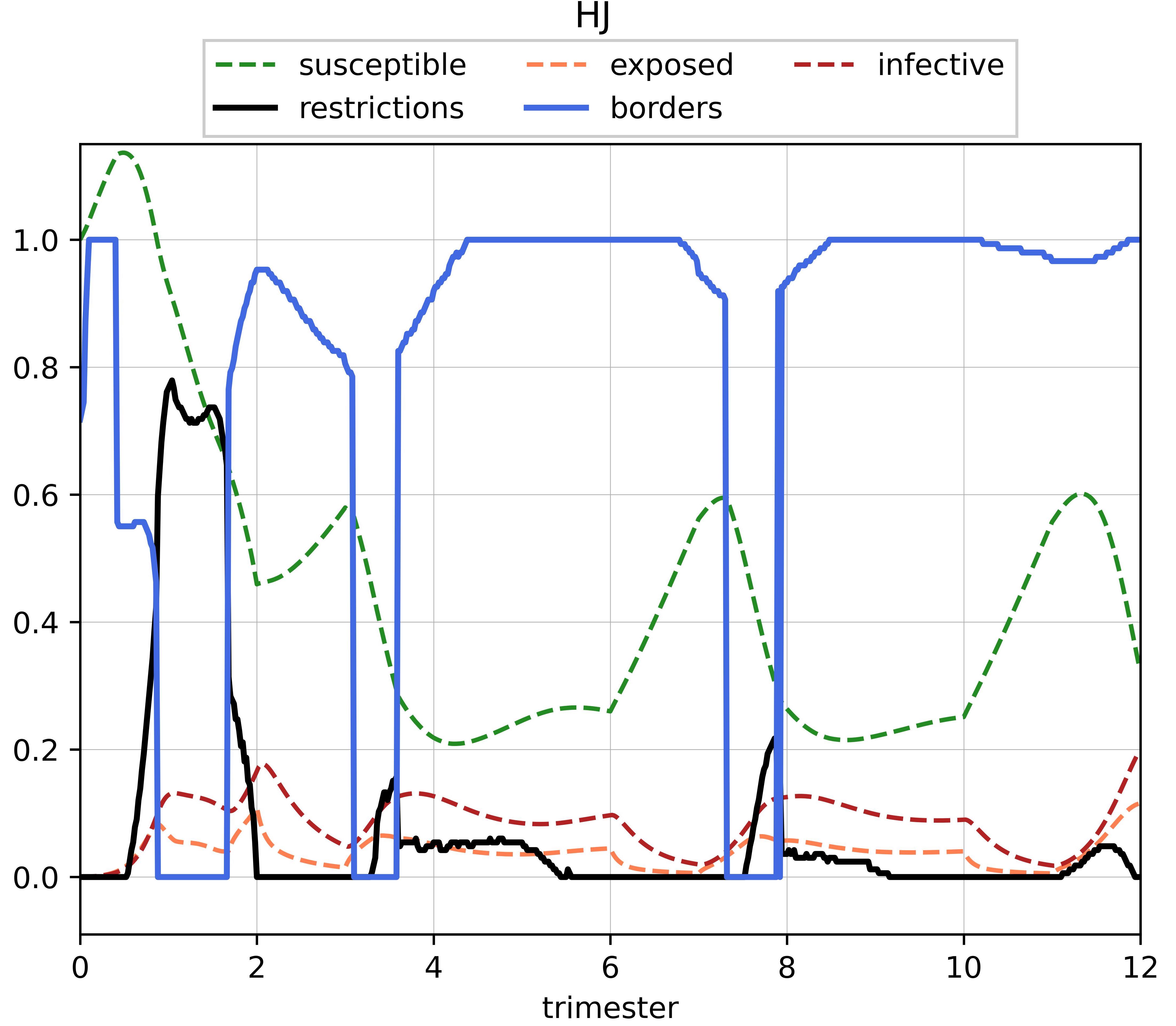}
    \includegraphics[width=7cm, height=6.1cm]{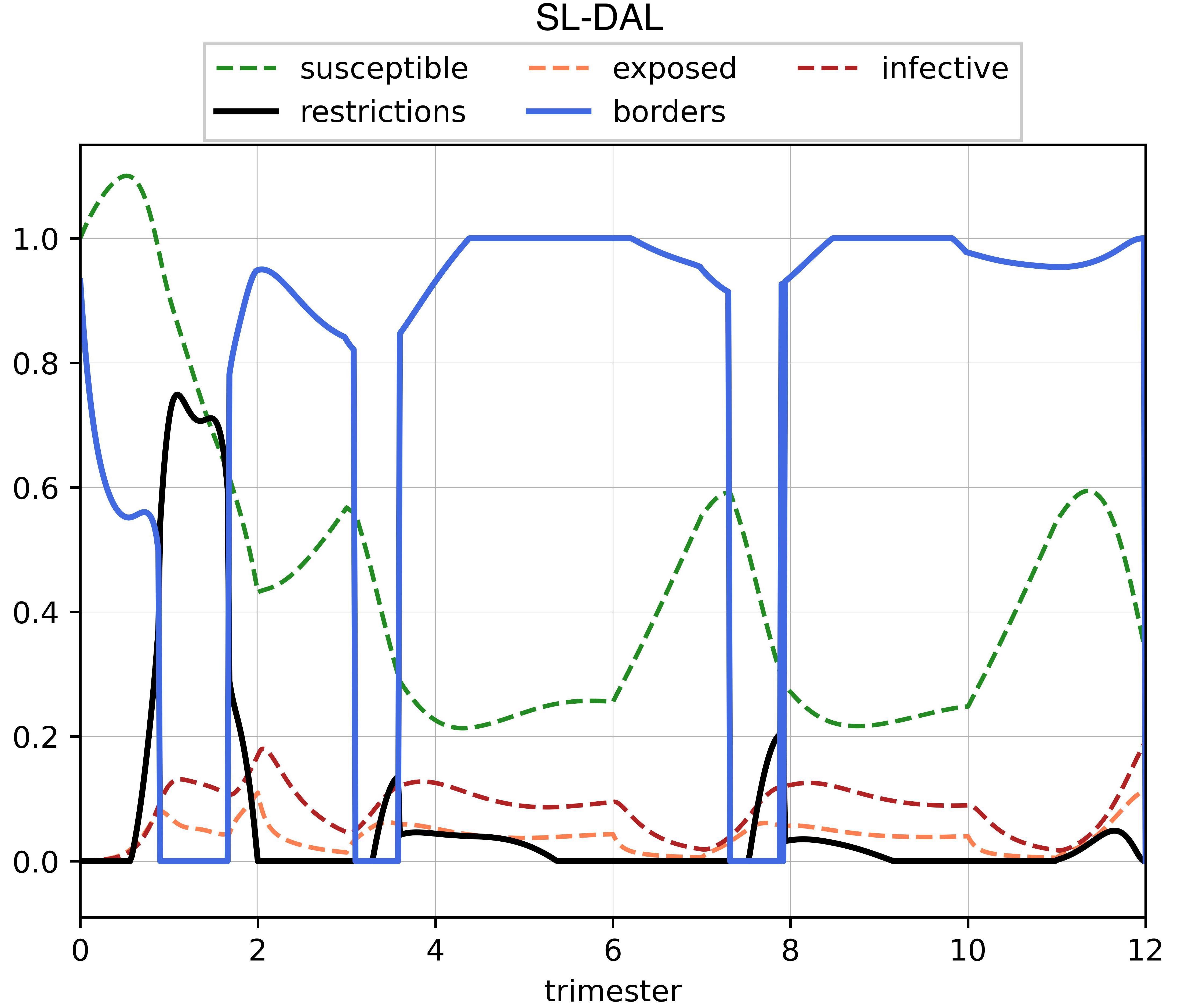}
    \caption{\textit{Test 3.} Optimal trajectories (dashed lines) and controls (full lines) for system \eqref{seir_frontiere}. Solutions obtained with the semi-Lagrangian scheme on the left, with Direct-Adjoint Looping, taking the solution given by the SL scheme as initial guess for the controls, on the right.}
    \label{fig:sim3}
\end{figure}

We perform a first test with the DAL algorithm, selecting as initial guess for the controls $\lambda^{(0)}=[0,\ldots,0]$ and $\phi^{(0)}=[1,\ldots,1]$. This choice corresponds to the uncontrolled scenario, where no restrictions are applied and the borders are open. The results are shown on the left part of Figure \ref{fig:sim3wrong} and we will indicate them with the superscript \textit{left}. We observe an initial lockdown, followed by three other periods of mild mobility restrictions; the borders are shut only around the first trimester and after that they are partially closed twice. We then repeat the test changing the initial guess on the borders to $\phi^{(0)}=[0,\ldots,0]$, obtaining the results reported on the right side of Figure \ref{fig:sim3wrong}, which we will indicate with the superscript \textit{right}. They are drastically different from those obtained in the first run and resemble the solution relative to the basic model in Figure \ref{fig:sim1}. In fact, the algorithm is keeping the borders closed all the time, and setting $\phi(\tau)=0$ for all $\tau$ is precisely equivalent to the basic model \eqref{seir_controllato_ridotto} with no vaccination. This means that at least one of the two solutions we obtained is only a local minimum (or just a stationary point) of the cost functional. A quick comparison between the value of the cost functional computed along the two solutions confirms that the first is cheaper than the second: $J_{x,t}^{\text{\textit{left}}}=20.092728$, whereas $J_{x,t}^{\text{\textit{right}}}=22.321380$. As a consequence, we can be sure that the solution obtained with $\phi^{(0)}=[0,\ldots,0]$ is not optimal. At this point, we run the semi-Lagrangian scheme, obtaining some discrete controls $\lambda^{SL}$ and $\phi^{SL}$, and then initialise the DAL algorithm with $\lambda^{(0)}=\lambda^{SL}$ and $\phi^{(0)}=\phi^{SL}$. In this way, we obtain the solutions reported in Figure \ref{fig:sim3}. We observe that the two solutions are rather different from the ones we had obtained before: there are three distinct periods in which the borders are shut, while the restrictions are still applied four times, more strictly in the beginning and then lowered progressively. The cost functional computed along the solutions given by the SL scheme and the combined SL-DAL algorithm is equal to, respectively, $J_{x,t}^{SL}=19.988674$ and $J_{x,t}^{SL-DAL}=19.977807$. This final check confirms that the actual minimum of the cost functional was not computed in any of the first two simulations, due to the different initialisation of the discrete controls, and that the combination with the semi-Lagrangian scheme can really make a difference. 

\subsection{Test 4: state constraints}

For this simulation we go back to the basic model, but we adopt a different approach to the problem. Instead of penalising the infective through their cost, we impose a state constraint on the system. This is not an abstract hypothesis, since every country has a limited number of hospital beds and, most importantly, of Intensive Care Units (ICU). We assume, for simplicity, that a fixed percentage $p_{ICU}$ of those who contract the disease end up in intensive care. This hypothesis can be misleading in models with heterogeneity, but since we are working with the epidemic SEIR model, the fundamental assumption is that the population is homogeneous and well mixed, as already mentioned in Section \ref{sec:SEIR}. For this reason, our assumption is reasonable and we set $p_{ICU}=0.05\%$. In this way, the constraint on intensive care units translates into a constraint on the state variable $i$. Taking the total number of intensive care units available in Italy before the Covid-19 pandemic \cite{terapie_intensive} and normalising it as we previously did for the initial data, we obtain an upper bound on $i(\tau)$:
$$
i_{\max} = 0.13.
$$

We keep only the components of the running cost \eqref{running_cost} related to restrictions and vaccines and add the penalisation term discussed in Section \ref{subsec:constraints}. Once again, we run the DAL algorithm with various initial guesses, always obtaining the trajectories and controls reported in Figure \ref{fig:sim4}.
\begin{figure}[h]
    \centering
    \includegraphics[width=7cm, height=7cm]{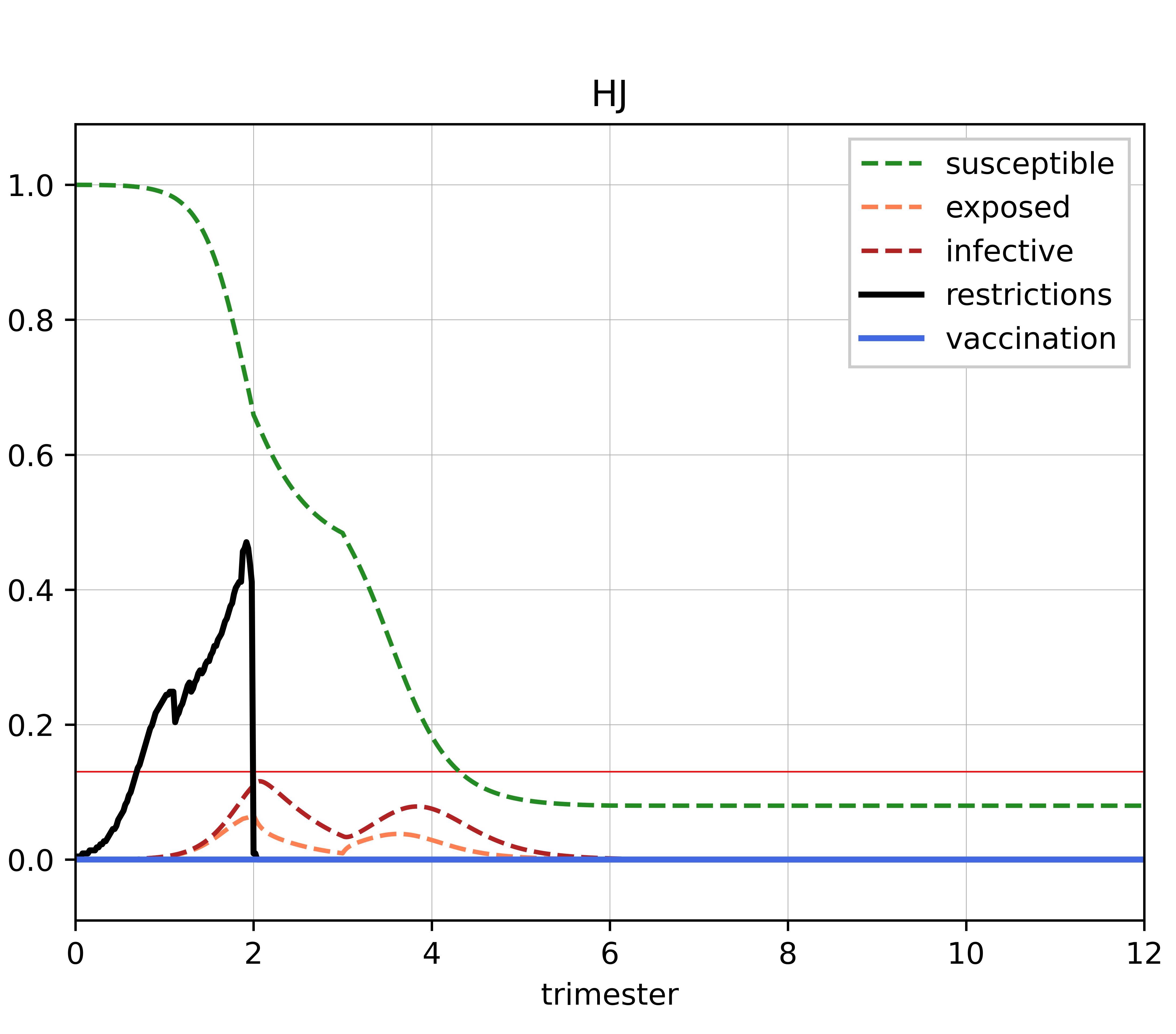}
    \includegraphics[width=7cm, height=7cm]{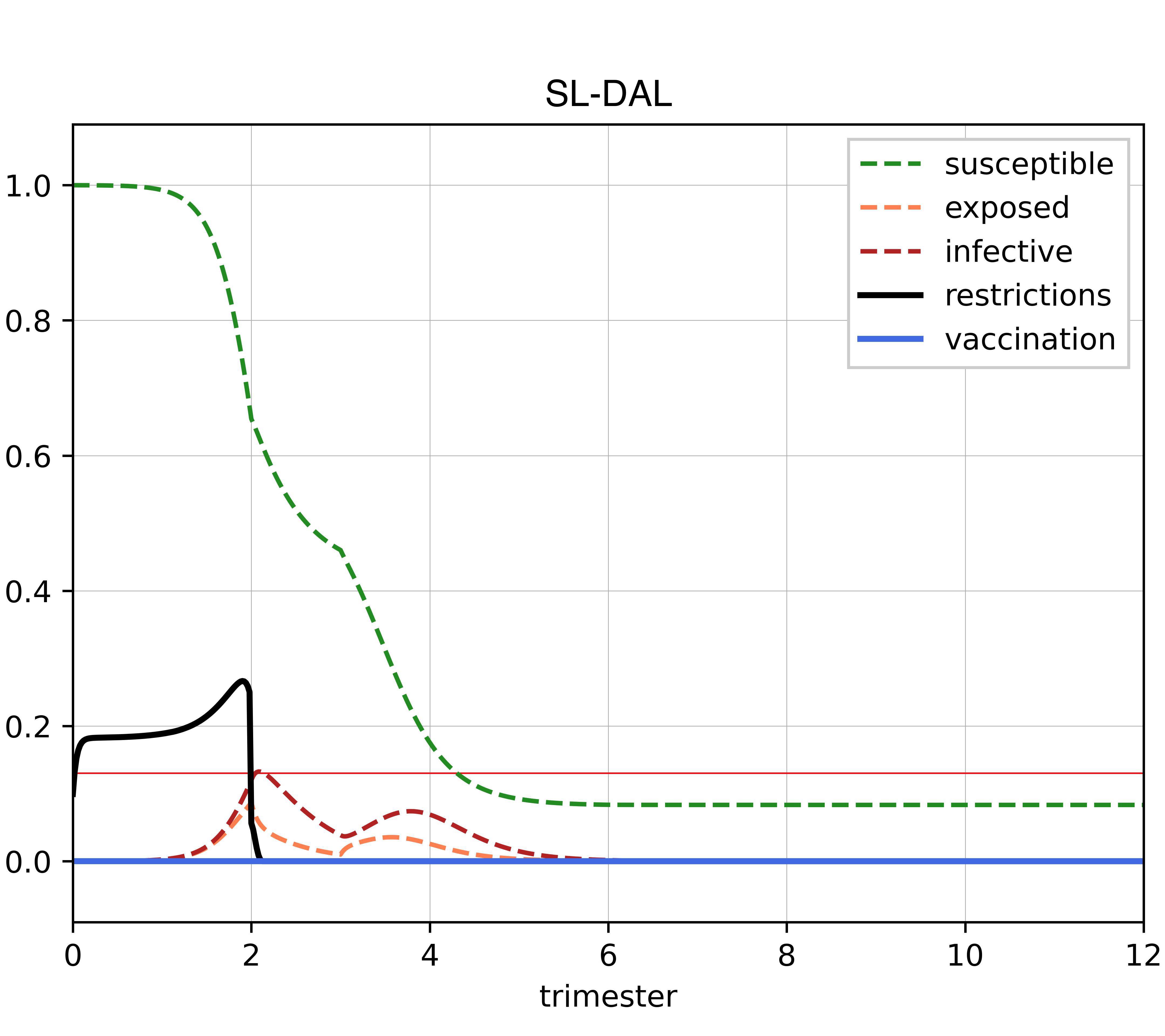}
    \caption{\textit{Test 4.} Optimal trajectories (dashed lines) and controls (full lines) for system \eqref{seir_controllato_ridotto} with a state constraint on $i$. Solutions obtained with the semi-Lagrangian scheme on the left, with Direct-Adjoint Looping on the right. The horizontal red line represents the threshold $i_{\max}$.}
    \label{fig:sim4}
\end{figure}
Comparing this scenario with that of \textit{Test 1}, without the state constraint, we notice that the strategy is not the exact same. Here we have a six-month period of mobility restrictions as soon as the disease starts spreading, but still no vaccination as in the unconstrained case. Since in this case the goal is not to minimise the number of infective at each time, but only to keep them strictly below the threshold $i_{\max}$, there is no need for further restrictions after the second trimester.

\subsection{Test 5: temporary immunity and state constraints}

Similarly to what we did for the previous simulation, we impose the same state constraint on $i$, representing the limited availability of Intensive Care Units, on the model with temporary immunity \eqref{seir_immunità_temporanea}. We use the same running cost of \textit{Test 4}, run the DAL algorithm with $\lambda^{(0)}=\nu^{(0)}=[0,\ldots,0]$ and obtain the solution reported in Figure \ref{fig:sim5wrong}.
\begin{figure}[!h]
    \centering
    \includegraphics[width=7cm]{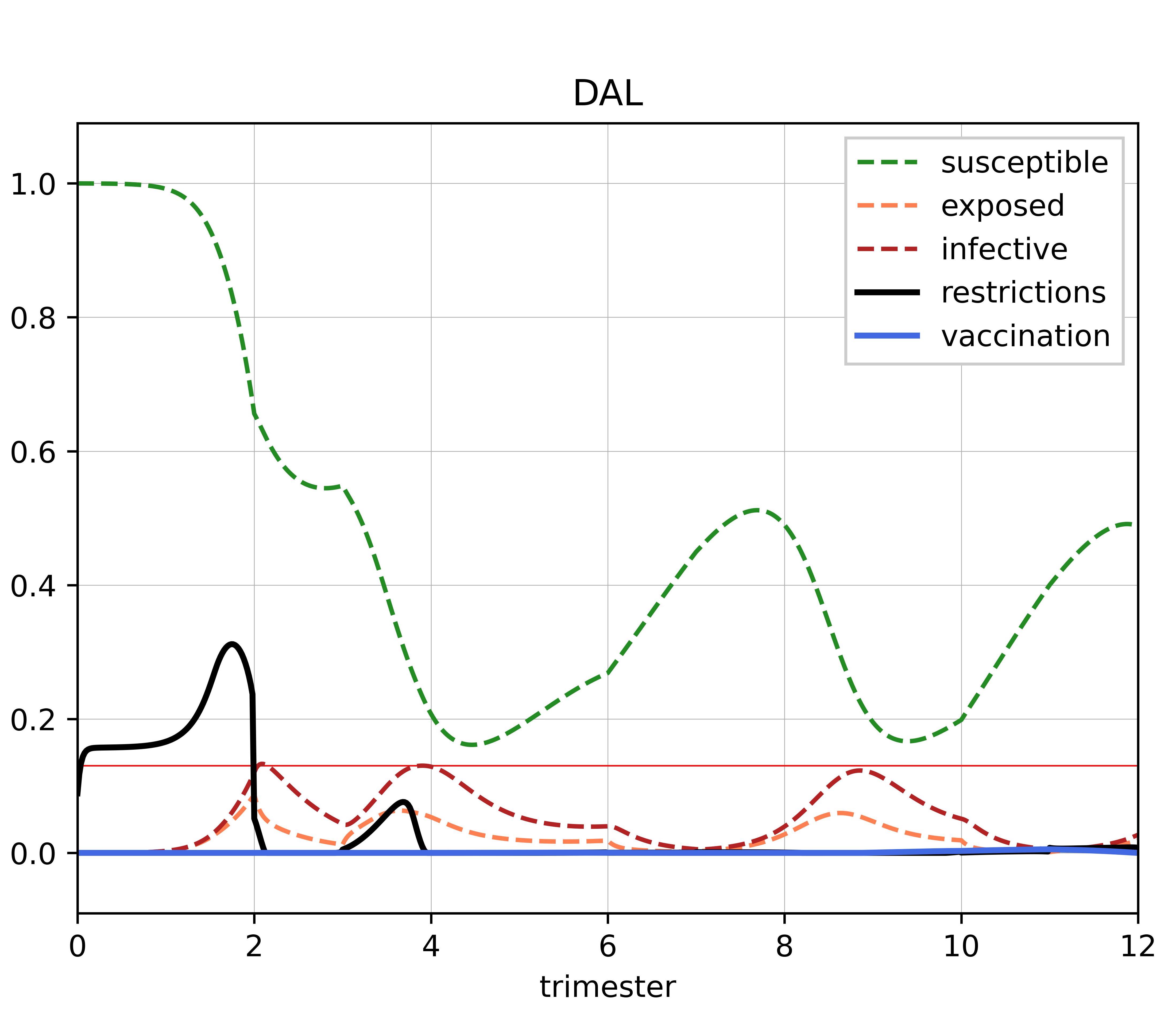}
    \caption{\textit{Test 5.} Optimal trajectories (dashed lines) and controls (full lines) for system \eqref{seir_immunità_temporanea} with a state constraint on $i$. Solutions obtained with the with the Direct-Adjoint Looping initialised with $\lambda^{(0)}=\nu^{(0)}=[0,\ldots,0]$. The horizontal line represents the threshold $i_{\max}$.}
    \label{fig:sim5wrong}
\end{figure}
We observe two restriction periods, one in the first six months and another one during the fourth trimester, while vaccines are practically never used. The cost functional evaluated along this solution is $J_{x,t}^{DAL}=0.362150$. We then run the SL scheme and the DAL algorithm initialised with $\lambda^{(0)}=\lambda^{SL}$ and $\nu^{(0)}=\nu^{SL}$, obtaining the results shown in Figure \ref{fig:sim5}.
\begin{figure}[!h]
    \centering
    \includegraphics[width=7cm, height=7cm]{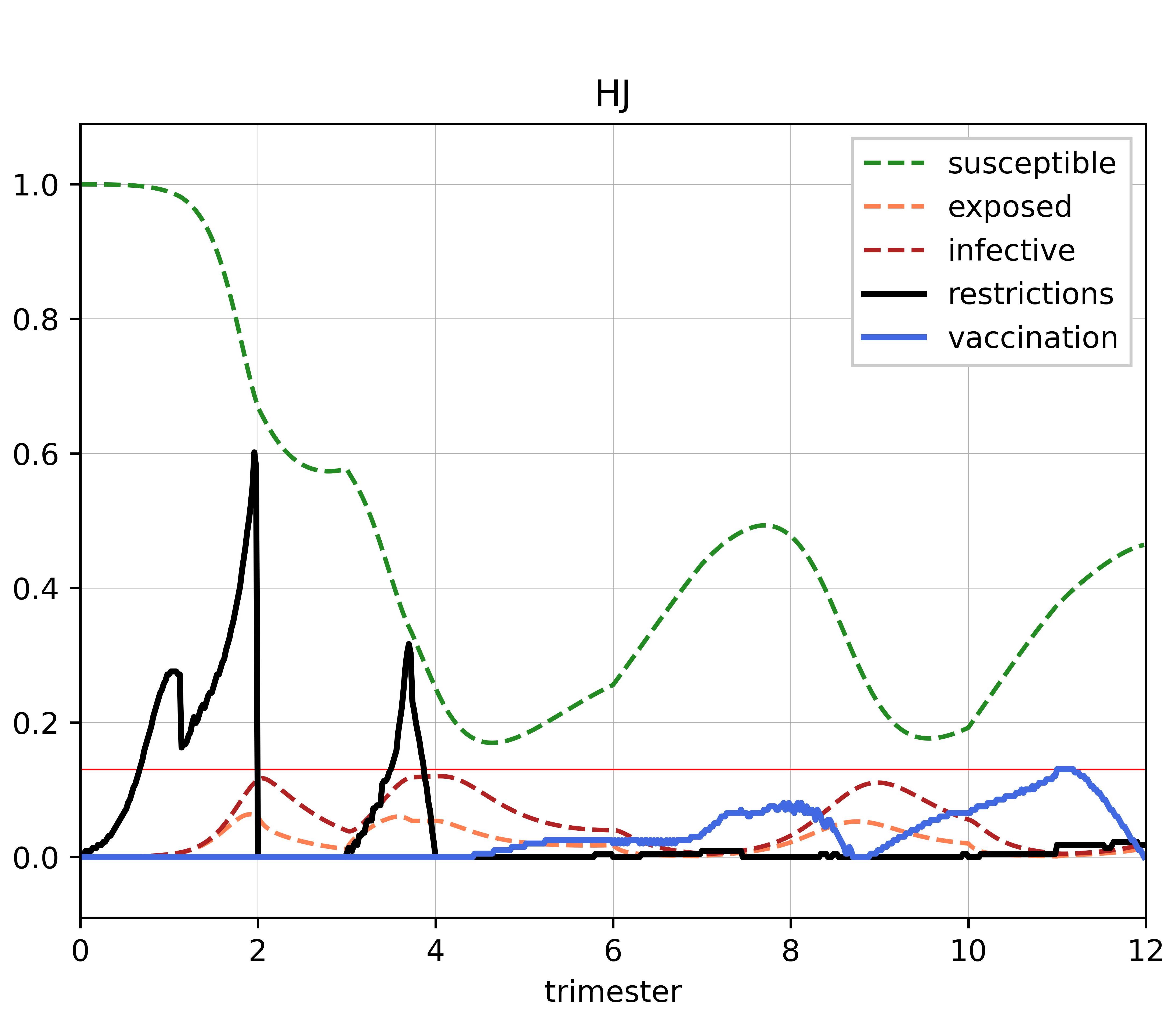}
    \includegraphics[width=7cm, height=7cm]{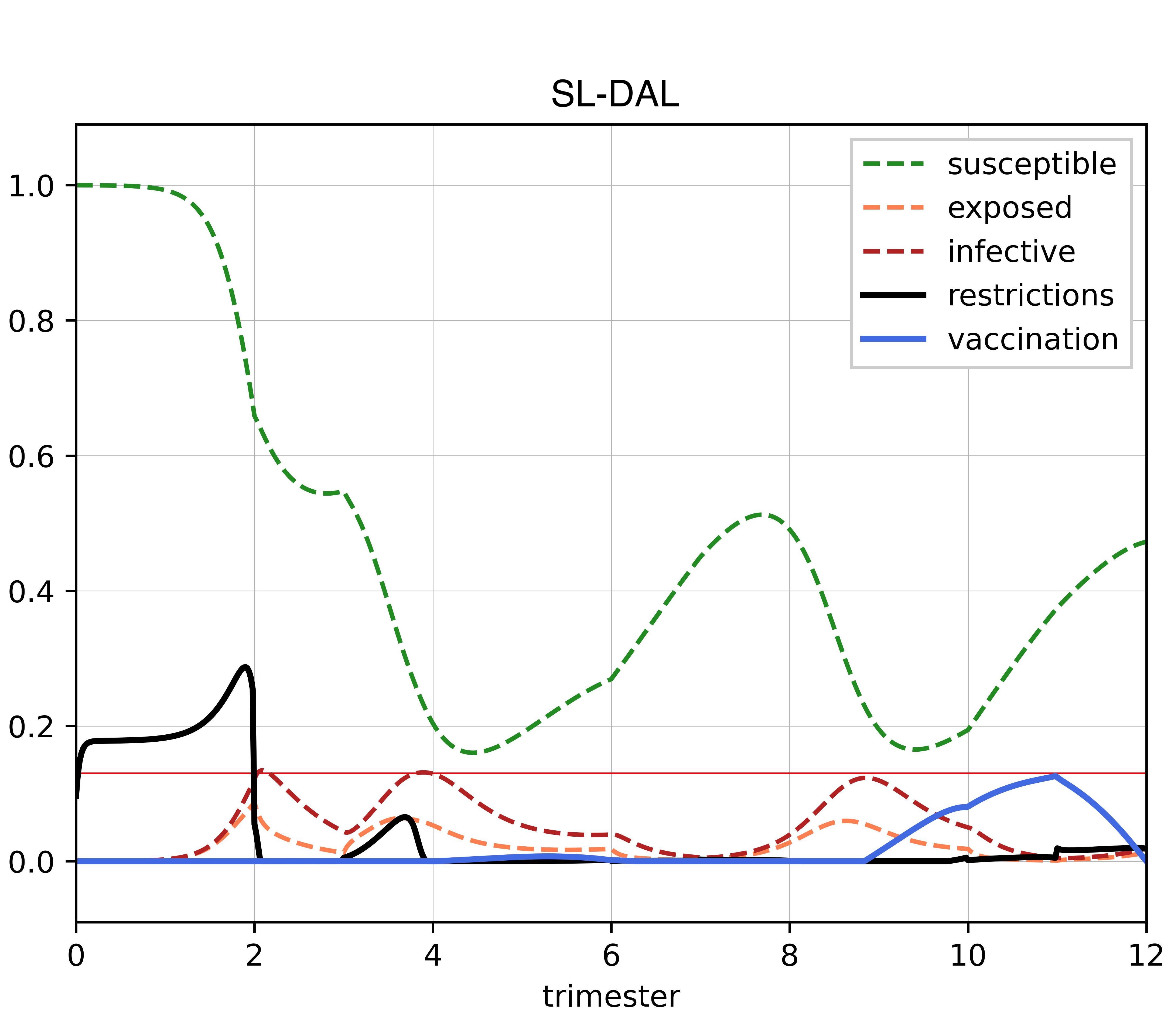}
    \caption{\textit{Test 5.} Optimal trajectories (dashed lines) and controls (full lines) for system \eqref{seir_immunità_temporanea} with a state constraint on $i$. Solutions obtained with the semi-Lagrangian scheme on the left, with Direct-Adjoint Looping initialised with $\lambda^{(0)}=\lambda^{SL}$ and $\nu^{(0)}=\nu^{SL}$ on the right. The horizontal line represents the threshold $i_{\max}$.}
    \label{fig:sim5}
\end{figure}
 Similarly to what we found in \textit{Test 3}, the solution we obtain is different from the one we get for $\lambda^{(0)}=\nu^{(0)}=[0,\ldots,0]$. Although the profile of $\lambda(\tau)$ is similar to what we obtained before, there is a remarkable difference in $\nu(\tau)$. The cost of this solution is $J_{x,t}^{SL-DAL}=0.333978 < J_{x,t}^{DAL}$, confirming that this is indeed the minimum of the cost functional.

\subsection{Summary}
For the sake of completeness, we summarise in Table \ref{tab:costi} the results of all the tests we performed. In particular, we report the values of the cost functional $J_{x,t}$ evaluated along, respectively, the trajectories of the uncontrolled system ($J_{x,t}^U$), those computed by the semi-Lagrangian scheme ($J_{x,t}^{SL}$), and those obtained with the Direct-Adjoint Looping method initialised with the output of the SL scheme ($J_{x,t}^{SL-DAL}$). In addition, we report the discrete $L^\infty$ norm of the difference between the optimal trajectories obtained with the two algorithms:
$$
|| \Delta y ||_\infty := 
\begin{pmatrix}
    || s^{SL} - s^{SL-DAL} ||_\infty \\
    || e^{SL} - e^{SL-DAL} ||_\infty \\
    || i^{SL} - i^{SL-DAL} ||_\infty
\end{pmatrix}\tpose.
$$
\begin{table}[h]
    \centering
    \begin{tabular}{|c|c|c|c|c|}
        \hline
        \rule{0pt}{10pt} & $J_{x,t}^U$ & $J_{x,t}^{SL}$ & $J_{x,t}^{SL-DAL}$ & $|| \Delta y ||_\infty$  \\
        \hline
        \textit{Test 1} & 20.990463 & 20.526586 & 20.521155 & (0.028566, 0.005591, 0.009907) \\
        \hline
        \textit{Test 2} & 20.180178 & 19.897859 & 19.865984 & (0.095190, 0.010398, 0.018655) \\
        \hline
        \textit{Test 3} & 21.443343 & 19.988674 & 19.977807 &  (0.039488, 0.008750, 0.013415) \\
        \hline
        \textit{Test 4} & 6.413498 & 0.412412 & 0.296018 &  (0.033831, 0.020147, 0.017608) \\
        \hline
        \textit{Test 5} & 6.936510 & 0.530465 & 0.333978 & (0.051467, 0.027062, 0.020124) \\
        \hline
    \end{tabular}
    \caption{Costs of the uncontrolled scenario and of the optimal policy for each simulation of Section \ref{sec:numerical}, obtained with the semi-Lagrangian scheme and the Direct-Adjoint Looping method initialised with the output of the SL scheme. The rightmost column contains the discrete $L^\infty$ norm of the difference between the optimal trajectories obtained with the two algorithms.}
    \label{tab:costi}
\end{table}
It is worth noting that $|| \Delta y ||_\infty$ is always at most of order $O(10^{-2})$, which is compatible with the magnitude of the discretisation steps in the SL scheme. 

Finally, to further validate the solutions found by the combined scheme, we check first and second order necessary optimality conditions \eqref{firstOrderConditions}--\eqref{secondOrderCondition} for all the tests we performed in this section. In Figures \ref{fig:firstOrder_1}--\ref{fig:firstOrder_5} we can see that \eqref{firstOrderConditions} is always verified. We point out that for \textit{Tests 1, 2, 4, 5} we do not report $\frac{\partial\mathcal{H}}{\partial\nu} (\tau)$ for $\tau \in [0,4)$, since $\nu(\tau) = 0$ for $\tau \in [0,4)$ in force of \eqref{numax}. Regarding second-order conditions, we observe that for \textit{Tests 1, 2, 4, 5} the Hessian matrix $\partial^2\mathcal{H} / \partial\alpha^2$, where $\alpha=(\lambda,\nu)$, is diagonal and its eigenvalues are non-negative, making \eqref{secondOrderCondition} immediately verified. For \textit{Test 3}, instead, we have
\begin{equation} \label{hessian}
    \frac{\partial^2\mathcal{H}}{\partial\alpha^2} =
    \begin{pmatrix}
        2\, c_\lambda\, (1+\delta\, \tau\, \phi) & 2\, c_\lambda\, \lambda\, \delta\, \tau \\
        2\, c_\lambda\, \lambda\, \delta\, \tau & 2\, c_\phi\, (1+ \delta\, \tau\, (3\,\phi -2))
    \end{pmatrix},
\end{equation}
which is symmetrical, therefore diagonalisable with real eigenvalues. Since \eqref{hessian} is a $2 \times 2$ matrix, in order to check \eqref{secondOrderCondition_a}, it is sufficient to verify that its trace and determinant are non-negative whenever the controls are not on the boundary of the control set. We can see in Figure \ref{fig:secondOrder} that
$$
\begin{aligned}
  \Tr \left[ \frac{\partial^2\mathcal{H}}{\partial\alpha^2} \right] > 0 \quad \text{ and } \quad \det \left[\frac{\partial^2\mathcal{H}}{\partial\alpha^2} \right] > 0 \quad &\text{ if } \quad (\lambda,\phi) \in (0, 0.9)\times(0, 1), \\
  \frac{\partial^2\mathcal{H}}{\partial\phi^2} > 0 \quad &\text{ if } \quad (\lambda = 0 \lor \lambda = 0.9) \land \phi \in (0,1), \\
  \frac{\partial^2\mathcal{H}}{\partial\lambda^2} > 0 \quad &\text{ if } \quad \lambda \in (0, 0.9) \land (\phi = 0 \lor \phi = 1), 
\end{aligned}
$$
meaning \eqref{secondOrderCondition} is verified for \textit{Test 3} as well and the inequalities are strictly respected.
\begin{figure}[!h]
    \centering
    \includegraphics[width=15cm]{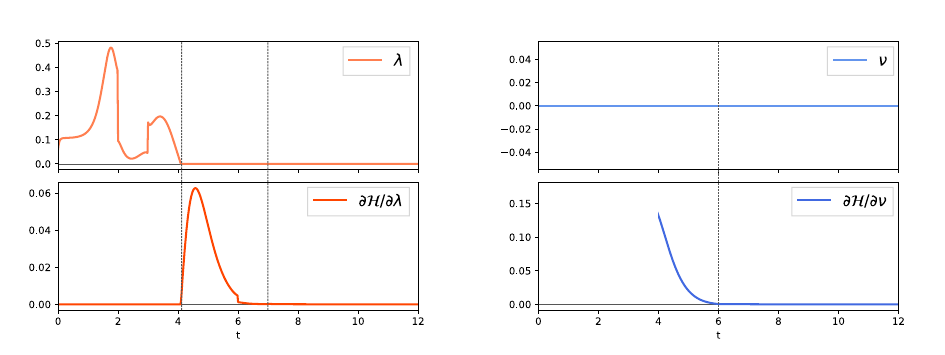}
    \caption{First-order necessary conditions \eqref{firstOrderConditions} for the solutions obtained with the combined SL-DAL approach for \textit{Test 1}.}
    \label{fig:firstOrder_1}
\end{figure}
\begin{figure}[!h]
    \centering
    \includegraphics[width=15cm]{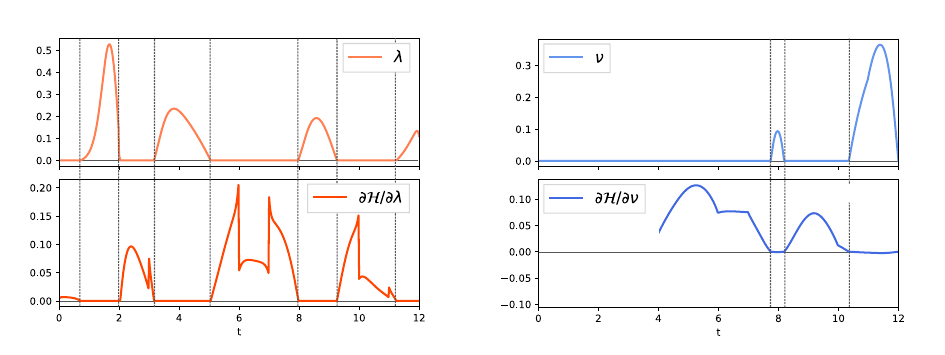}
    \caption{First-order necessary conditions \eqref{firstOrderConditions} for the solutions obtained with the combined SL-DAL approach for \textit{Test 2}.}
    \label{fig:firstOrder_2}
\end{figure}
\begin{figure}[!h]
    \centering
    \includegraphics[width=15cm]{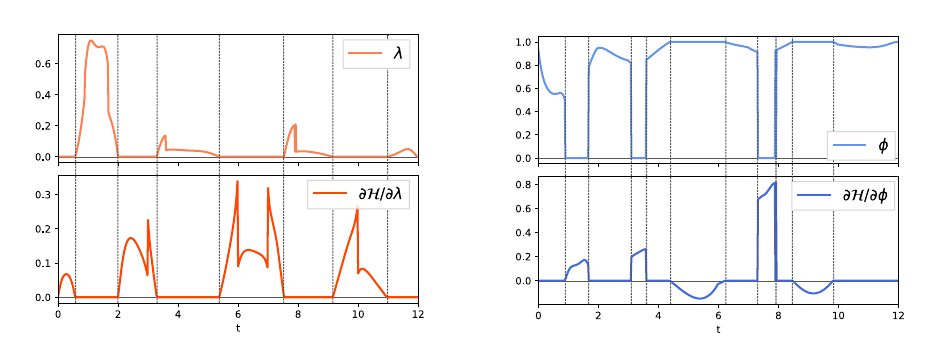}
    \caption{First-order necessary conditions \eqref{firstOrderConditions} for the solutions obtained with the combined SL-DAL approach for \textit{Test 3}.}
    \label{fig:firstOrder_3}
\end{figure}
\begin{figure}[!h]
    \centering
    \includegraphics[width=15cm]{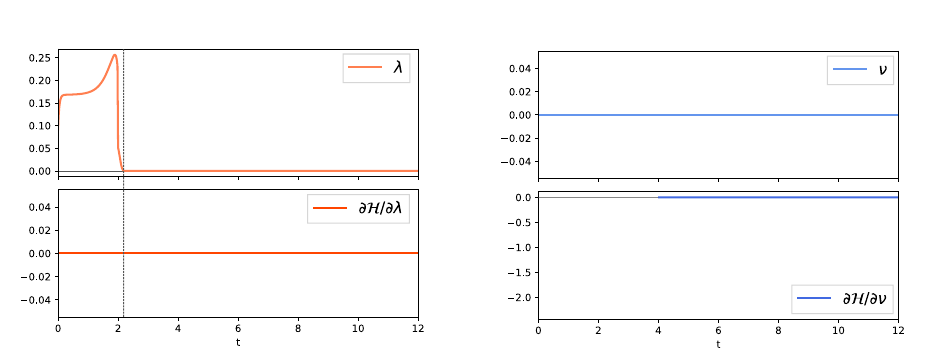}
    \caption{First-order necessary conditions \eqref{firstOrderConditions} for the solutions obtained with the combined SL-DAL approach for \textit{Test 4}.}
    \label{fig:firstOrder_4}
\end{figure}
\begin{figure}[!h]
    \centering
    \includegraphics[width=15cm]{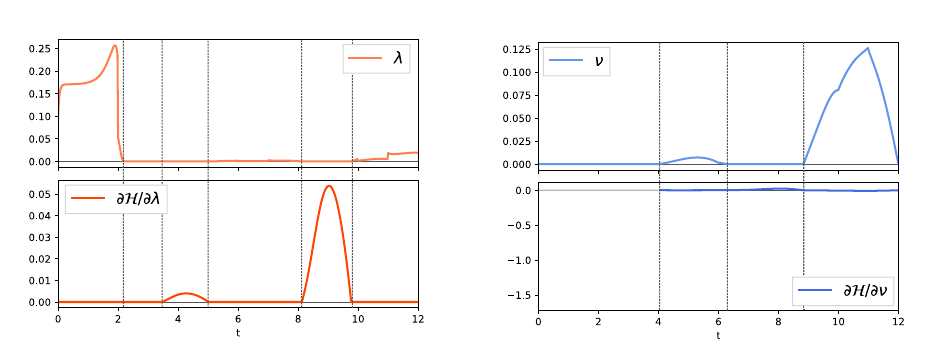}
    \caption{First-order necessary conditions \eqref{firstOrderConditions} for the solutions obtained with the combined SL-DAL approach for \textit{Test 5}.}
    \label{fig:firstOrder_5}
\end{figure}
\begin{figure}[!h]
    \centering
    \includegraphics[width=7cm]{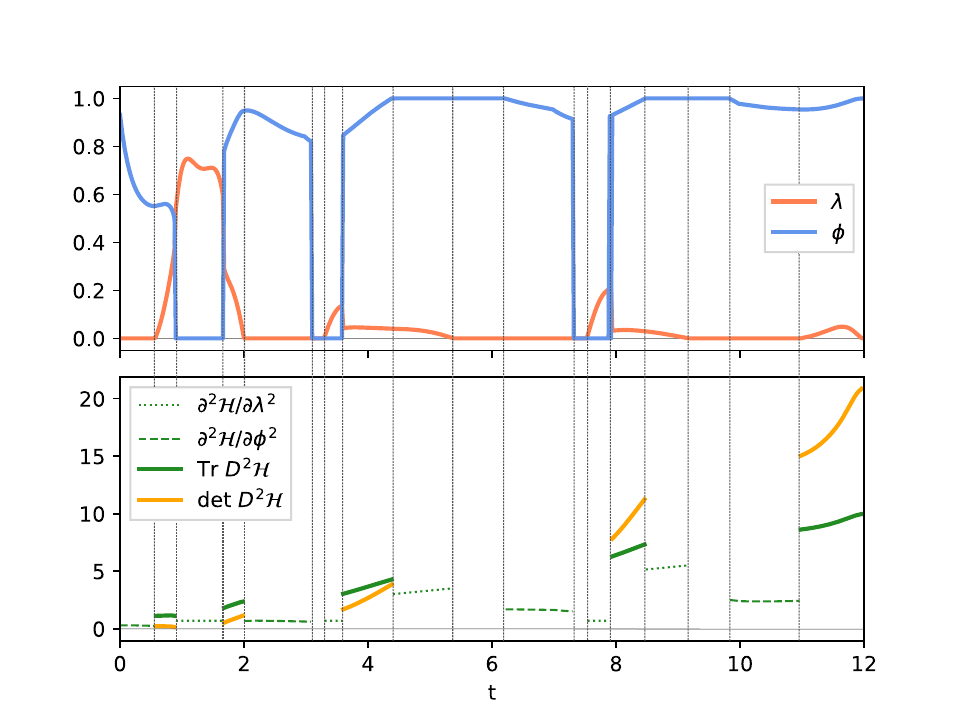}
    \caption{Second-order necessary conditions \eqref{secondOrderCondition} for the solutions obtained with the combined SL-DAL approach for \textit{Test 3}.}
    \label{fig:secondOrder}
\end{figure}\\

It is worth noting that, unlike the finite-dimensional setting, verifying \eqref{secondOrderCondition} strictly is still not sufficient to prove that the solution of \eqref{pontryagin_var_zero} is a local minimum of \eqref{hamiltonian}. In general, point-wise conditions such as the strong Legendre or Legendre-Clebsch condition, must always be coupled with additional ones, like the strong Jacobi condition or the existence of a bounded solution of a particular Riccati equation (see \cite{M_1981,CZ_1986,IT_2009,CT_2015} for reference). This, however, goes beyond the scope of this paper and will be addressed in future works.

The results above confirm that the Direct-Adjoint Looping algorithm, if initialised with the output of a semi-Lagrangian scheme, produces solutions that are in agreement with those given by the SL scheme and verify first and second order necessary optimality conditions. Moreover, as expected, $J_{x,t}^{SL-DAL} < J_{x,t}^{SL}$ for all simulations, since the state-space and the control set are not discretised in the variational approach, leading to more precise solutions. In particular, \textit{Test 3} and \textit{Test 5} are a concrete example in which the Dynamic Programming approach can significantly help to retrieve reliable control policies when the cost functional presents multiple stationary points.


\section{Conclusions}
\label{sec:conclusions}

We presented some variations of the epidemic SEIR model, including a variable transmission rate, the initial unavailability of vaccines, temporary immunity, state constraints on Intensive Care Units and interactions with external populations. We formulated some finite horizon optimal control problems associated with them, considering vaccines, restrictive measures and the possibility to close the borders as controls. We presented two different theoretical approaches, Dynamic Programming and Pontryagin's Maximum Principle, in order to devise suitable approximation procedures for the computation of optimal strategies. We performed several numerical simulations and showed that descent methods based on the variational approach can be highly sensitive to the initial guess on the controls, and this can lead to sub-optimal solutions. We showed that a combination of the two methods, where we initialise the descent algorithm with the solution given by the semi-Lagrangian scheme, can help to obtain high-quality, reliable approximations of the optimal controls. In order to test the quality of the solutions obtained by the combined SL-DAL algorithm, we checked numerically first and second order necessary optimality conditions.

The aim of this work is to give an idea of the potential and the advantages of the combination of the two approaches for optimal control problems in epidemiology. This idea can also be applied to more complex models, aiming at capturing more realistic scenarios. In this case, a collaboration with epidemiologists would be needed in order to estimate the parameters from real data, not only for the dynamical system, but for the cost functional as well. The same framework can be applied to other compartmental models and different controls can be considered, based on the particular characteristics of the infectious disease or on the containment instruments that are available in a certain area. Even more in general, the techniques we presented can be adapted to other differential systems in biomathematics for which it may be interesting to study some optimal control problems.

\section*{Declarations of interest}
None.

\section*{Acknowledgements}
The authors are members of the \textit{Gruppo Nazionale Calcolo Scientifico - Istituto Nazionale di Alta Matematica} (GNCS-INdAM).

This research was supported by the PNRR-MUR project \textit{Italian Research Center on High Performance Computing, Big Data and Quantum Computing}, by the Italian PNRR fund within the doctoral project \textit{Modelli matematici per la simulazione e controllo delle pandemie} at Sapienza University of Rome, and by the \textit{INdAM-GNCS Project} code CUP\_E53C23001670001.

The authors would like to dedicate this work to the memory of professor Maurizio Falcone, who passed away in November 2022. His long experience in numerical analysis and optimal control, together with his ideas and enthusiasm for research, were of great inspiration to start this research project.

\medskip

\bibliographystyle{siam} 
\bibliography{seir.bib}

\begin{thebibliography}{10}

\bibitem{terapie_intensive}
{\sc {Agenzia Nazionale per i Servizi Sanitari Regionali}}.
\newblock
  \url{https://www.agenas.gov.it/covid19/web/index.php?r=site}\texttt{$\%$2Fgraph3},
  2022.

\bibitem{AFS_2020}
{\sc A.~Alla, M.~Falcone, and L.~Saluzzi}, {\em A tree structure algorithm for
  optimal control problems with state constraints}, Rendiconti di Matematica e
  delle sue Applicazioni, 41 (2020), pp.~193--221.

\bibitem{AP_2022}
{\sc R.~Arbel and J.~Pliskin}, {\em Vaccinations versus lockdowns to prevent
  {COVID-19} mortality}, Vaccines, 10 (2022).

\bibitem{AS_1984}
{\sc J.~L. Aron and I.~B. Schwartz}, {\em Seasonality and period-doubling
  bifurcations in an epidemic model.}, Journal of Theoretical Biology, 110
  (1984), pp.~665--679.

\bibitem{AFG_2022}
{\sc F.~Avram, L.~Freddi, and D.~Goreac}, {\em Optimal control of a {SIR}
  epidemic with {ICU} constraints and target objectives}, Applied Mathematics
  and Computation, 418 (2022), p.~126816.

\bibitem{AKK_2021}
{\sc B.~Azmi, D.~Kalise, and K.~Kunisch}, {\em Optimal feedback law recovery by
  gradient-augmented sparse polynomial regression}, Journal of Machine Learning
  Research, 22 (2021), pp.~1--32.

\bibitem{B_2011}
{\sc N.~Baca{\"e}r}, {\em A Short History of Mathematical Population Dynamics},
  Springer London, 2011.

\bibitem{BCD_2008}
{\sc M.~Bardi and I.~Capuzzo-Dolcetta}, {\em Optimal Control and Viscosity
  Solutions of {Hamilton}-{Jacobi}-{Bellman} Equations}, Modern Birkh{\"a}user
  Classics, Birkh{\"a}user Boston, 2008.

\bibitem{Beh_2000}
{\sc H.~Behncke}, {\em Optimal control of deterministic epidemics}, Optimal
  Control Applications and Methods, 21 (2000), pp.~269--285.

\bibitem{B_1957}
{\sc R.~Bellman}, {\em {Dynamic Programming}}, Princeton University Press,
  Princeton, NJ, 1957.

\bibitem{BDPV_2021}
{\sc P.-A. Bliman, M.~Duprez, Y.~Privat, and N.~Vauchelet}, {\em Optimal
  immunity control and final size minimization by social distancing for the
  {SIR} epidemic model}, Journal of Optimization Theory and Applications, 189
  (2021), pp.~408--436.

\bibitem{BGY_1997}
{\sc S.~Blount, A.~Galambosi, and S.~Yakowitz}, {\em Nonlinear and dynamic
  programming for epidemic intervention}, Applied Mathematics and Computation,
  86 (1997), pp.~123--136.

\bibitem{BMO_2019}
{\sc B.~Buonomo, P.~Manfredi, and A.~d’Onofrio}, {\em Optimal time-profiles
  of public health intervention to shape voluntary vaccination for childhood
  diseases}, Journal of mathematical biology, 78 (2019), pp.~1089--1113.

\bibitem{CCFP_2012}
{\sc S.~Cacace, E.~Cristiani, M.~Falcone, and A.~Picarelli}, {\em A patchy
  dynamic programming scheme for a class of {Hamilton-Jacobi-Bellman}
  equations}, {SIAM} Journal on Scientific Computing, 34 (2012),
  pp.~A2625--A2649.

\bibitem{CF_2021}
{\sc S.~Cacace and R.~Ferretti}, {\em Efficient implementation of
  characteristic-based schemes on unstructured triangular grids}, Computational
  and Applied Mathematics, 41 (2021).

\bibitem{CLL_2023}
{\sc S.~Cacace, A.~C. Lai, and P.~Loreti}, {\em A dynamic programming approach
  for controlled fractional {SIS} models}, Nonlinear Differential Equations and
  Applications NoDEA, 30 (2023), p.~20.

\bibitem{CT_2015}
{\sc E.~Casas and F.~Tr{\"o}ltzsch}, {\em Second order optimality conditions
  and their role in pde control}, Jahresbericht der Deutschen
  Mathematiker-Vereinigung, 117 (2015), pp.~3--44.

\bibitem{FCAAHT_2023}
{\sc C.~Chukwu, R.~Alqahtani, C.~Alfiniyah, F.~Herdicho, et~al.}, {\em A
  {Pontryagin’s} maximum principle and optimal control model with
  cost-effectiveness analysis of the {COVID-19} epidemic}, Decision Analytics
  Journal, 8 (2023), p.~100273.

\bibitem{CZ_1986}
{\sc F.~H. Clarke and V.~Zeidan}, {\em Sufficiency and the jacobi condition in
  the calculus of variations}, Canadian Journal of Mathematics, 38 (1986),
  pp.~1199--1209.

\bibitem{CEL_1984}
{\sc M.~G. Crandall, L.~C. Evans, and P.-L. Lions}, {\em Some properties of
  viscosity solutions of {Hamilton-Jacobi} equations}, Transactions of the
  American Mathematical Society, 282 (1984), pp.~487--502.

\bibitem{CL_1983}
{\sc M.~G. Crandall and P.-L. Lions}, {\em Viscosity solutions of
  {Hamilton-Jacobi} equations}, Transactions of the American mathematical
  society, 277 (1983), pp.~1--42.

\bibitem{CM_2010}
{\sc E.~Cristiani and P.~Martinon}, {\em Initialization of the shooting method
  via the {Hamilton-Jacobi-Bellman} approach}, Journal of Optimization Theory
  and Applications, 146 (2010), pp.~321--346.

\bibitem{DO_2021}
{\sc R.~Della~Marca and A.~d’Onofrio}, {\em Volatile opinions and optimal
  control of vaccine awareness campaigns: Chaotic behaviour of the
  forward-backward sweep algorithm vs. heuristic direct optimization},
  Communications in Nonlinear Science and Numerical Simulation, 98 (2021),
  p.~105768.

\bibitem{DHB_2013}
{\sc O.~Diekmann, H.~Heesterbeek, and T.~Britton}, {\em Mathematical Tools for
  Understanding Infectious Disease Dynamics}, Princeton University Press, 2013.

\bibitem{DH_2000}
{\sc O.~Diekmann and J.~Heesterbeek}, {\em Mathematical epidemiology of
  infectious diseases: Model building, analysis and interpretation}, Wiley
  Series in Mathematical and Computational Biology, Chichester, Wiley,  (2000).

\bibitem{OKS_2022}
{\sc F.~El~Ouardighi, E.~Khmelnitsky, and S.~P. Sethi}, {\em Epidemic control
  with endogenous treatment capability under popular discontent and social
  fatigue}, Production and Operations Management, 31 (2022), pp.~1734--1752.

\bibitem{E_1998}
{\sc L.~Evans}, {\em Partial Differential Equations}, Graduate studies in
  mathematics, American Mathematical Society, 1998.

\bibitem{FF_SIAM2013}
{\sc M.~Falcone and R.~Ferretti}, {\em Semi-Lagrangian Approximation Schemes
  for Linear and {Hamilton}-{Jacobi} Equations}, Society for Industrial and
  Applied Mathematics, Philadelphia, PA, 2013.

\bibitem{FF_2016}
{\sc M.~Falcone and R.~Ferretti}, {\em Numerical methods for
  {Hamilton}-{Jacobi} type equations}, in Handbook of Numerical Methods for
  Hyperbolic Problems, R.~Abgrall and C.-W. Shu, eds., vol.~17 of Handbook of
  Numerical Analysis, Elsevier, 2016, pp.~603--626.

\bibitem{FG_1999}
{\sc M.~Falcone and T.~Giorgi}, {\em An approximation scheme for evolutive
  {Hamilton-Jacobi} equations}, in Stochastic Analysis, Control, Optimization
  and Applications: A Volume in Honor of W.H. Fleming, W.~M. McEneaney, G.~G.
  Yin, and Q.~Zhang, eds., Birkh{\"a}user Boston, Boston, MA, 1999,
  pp.~289--303.

\bibitem{FKS_2023}
{\sc M.~Falcone, G.~Kirsten, and L.~Saluzzi}, {\em Approximation of optimal
  control problems for the {Navier-Stokes} equation via multilinear {HJB-POD}},
  Applied Mathematics and Computation, 442 (2023), p.~127722.

\bibitem{FSMTMPFB_2022}
{\sc R.~R.~A. Fernandes, M.~da~Silva~Santos, C.~A. da~Silva~Magliano, B.~R.
  Tura, L.~S. D.~N. Macedo, M.~P. Padila, A.~C.~W. França, and A.~A. Braga},
  {\em Cost utility of vaccination against {COVID-19} in brazil}, Value in
  Health Regional Issues, 31 (2022), pp.~18--24.

\bibitem{FR_1975}
{\sc W.~Fleming and R.~Rishel}, {\em Deterministic and Stochastic Optimal
  Control}, Springer New York, 1975.

\bibitem{FLZ_2009}
{\sc S.~Fomel, S.~Luo, and H.~Zhao}, {\em Fast sweeping method for the factored
  eikonal equation}, Journal of Computational Physics, 228 (2009),
  pp.~6440--6455.

\bibitem{GDA_2020}
{\sc G.~Gonzalez-Parra, M.~D{\'\i}az-Rodr{\'\i}guez, and A.~J. Arenas}, {\em
  Mathematical modeling to design public health policies for {Chikungunya}
  epidemic using optimal control}, Optimal Control Applications and Methods, 41
  (2020), pp.~1584--1603.

\bibitem{GS_2014}
{\sc B.-Z. Guo and B.~Sun}, {\em Dynamic programming approach to the numerical
  solution of optimal control with paradigm by a mathematical model for drug
  therapies of {HIV/AIDS}}, Optimization and Engineering, 15 (2014),
  pp.~119--136.

\bibitem{HD_2011}
{\sc E.~Hansen and T.~Day}, {\em Optimal control of epidemics with limited
  resources}, Journal of Mathematical Biology, 62 (2011), pp.~423--451.

\bibitem{multipliers-H}
{\sc M.~Hestenes}, {\em Multiplier and gradient methods}, Journal of
  Optimization Theory and Applications, 4 (1969), pp.~303--320.

\bibitem{HWL_2020}
{\sc Y.~G. Hwang, H.-D. Kwon, and J.~Lee}, {\em Feedback control problem of an
  {SIR} epidemic model based on the {Hamilton-Jacobi-Bellman} equation},
  Mathematical Biosciences and Engineering, 17 (2020), pp.~2284--2301.

\bibitem{IS_2013}
{\sc D.~Iacoviello and N.~Stasio}, {\em Optimal control for {SIRC} epidemic
  outbreak}, Computer methods and programs in biomedicine, 110 (2013),
  pp.~333--342.

\bibitem{IT_2009}
{\sc A.~D. Ioffe and V.~M. Tihomirov}, {\em Theory of Extremal Problems},
  Elsevier, 2009.

\bibitem{ISTAT}
{\sc {Istituto Nazionale di Statistica}}.
\newblock \url{http://dati.istat.it/Index.aspx?QueryId=18460}, 2022.

\bibitem{KS_1987}
{\sc C.~Kelley and E.~Sachs}, {\em Quasi-{Newton} methods and unconstrained
  optimal control problems}, {SIAM} Journal on Control and Optimization, 25
  (1987), pp.~1503--1516.

\bibitem{KMcK_1927}
{\sc W.~O. Kermack and A.~G. McKendrick}, {\em A contribution to the
  mathematical theory of epidemics}, Proceedings of the Royal Society of
  London. Series A, Containing Papers of a Mathematical and Physical Character,
  115 (1927), pp.~700--721.

\bibitem{K_2004}
{\sc D.~Kirk}, {\em Optimal Control Theory: An Introduction}, Dover Books on
  Electrical Engineering Series, Dover Publications, 2004.

\bibitem{LT_2015}
{\sc L.~Laguzet and G.~Turinici}, {\em Global optimal vaccination in the {SIR}
  model: properties of the value function and application to cost-effectiveness
  analysis}, Mathematical Biosciences, 263 (2015), pp.~180--197.

\bibitem{L_1979}
{\sc C.~Lefevre}, {\em Optimal control of the simple stochastic epidemic with
  variable recovery rates}, Mathematical Biosciences, 44 (1979), pp.~209--219.

\bibitem{LHLZWHY_2021}
{\sc X.~Liu, J.~Huang, C.~Li, Y.~Zhao, D.~Wang, Z.~Huang, and K.~Yang}, {\em
  The role of seasonality in the spread of {COVID-19} pandemic}, Environmental
  Research, 195 (2021), p.~110874.

\bibitem{LY_1973}
{\sc W.~P. London and J.~A. Yorke}, {\em Recurrent outbreaks of measles,
  chickenpox and mumps: I. seasonal variation in contact rates}, American
  Journal of Epidemiology, 98 (1973), pp.~453--468.

\bibitem{LS_2014}
{\sc A.~C. Lowen and J.~Steel}, {\em Roles of humidity and temperature in
  shaping influenza seasonality}, Journal of Virology, 88 (2014),
  pp.~7692--7695.

\bibitem{M_1981}
{\sc H.~Maurer}, {\em First and second order sufficient optimality conditions
  in mathematical programming and optimal control}, Springer, 1981.

\bibitem{NW_2006}
{\sc J.~Nocedal and S.~J. Wright}, {\em Numerical Optimization}, Springer, New
  York, NY, USA, 2e~ed., 2006.

\bibitem{PBGM_1964}
{\sc L.~S. Pontryagin, V.~G. Boltyanskii, R.~V. Gamkrelidze, and M.~E. F.},
  {\em The mathematical theory of optimal processes}, A Pergamon Press Book,
  The Macmillan Co., New York, 1964.

\bibitem{multipliers-P}
{\sc M.~Powell}, {\em A method for nonlinear constraints in minimization
  problems}, Academic Press, New York, NY, 1969.

\bibitem{AFS_2022}
{\sc L.~Saluzzi, A.~Alla, and M.~Falcone}, {\em Error estimates for a tree
  structure algorithm solving finite horizon control problems}, ESAIM: COCV, 28
  (2022), p.~69.

\bibitem{S_1999}
{\sc J.~A. Sethian}, {\em Fast marching methods}, SIAM Review, 41 (1999),
  pp.~199--235.

\bibitem{SH_2019}
{\sc M.~K.~S. Sihombing and H.~Hartono}, {\em Optimal control of the spread of
  {Dengue} fever using dynamic programming}, in 2019 International Conference
  on Information and Communications Technology (ICOIACT), IEEE, 2019,
  pp.~698--703.

\bibitem{STW_2019}
{\sc B.~Sun, Z.-Z. Tao, and Y.-Y. Wang}, {\em Dynamic programming viscosity
  solution approach and its applications to optimal control problems},
  Mathematics Applied to Engineering, Modelling, and Social Issues,  (2019),
  pp.~363--420.

\bibitem{WFCJHCLY_2021}
{\sc W.-C. Wang, J.~C.-Y. Fann, R.-E. Chang, Y.-C. Jeng, C.-Y. Hsu, H.-H. Chen,
  J.-T. Liu, and A.~M.-F. Yen}, {\em Economic evaluation for mass vaccination
  against {COVID-19}}, Journal of the Formosan Medical Association, 120 (2021),
  pp.~S95--S105.

\bibitem{XXXY_2017}
{\sc D.~Xu, X.~Xu, Y.~Xie, and C.~Yang}, {\em Optimal control of an {SIVRS}
  epidemic spreading model with virus variation based on complex networks},
  Communications in Nonlinear Science and Numerical Simulation, 48 (2017),
  pp.~200--210.

\bibitem{YD_2021}
{\sc I.~Yegorov and P.~Dower}, {\em Perspectives on characteristics based
  curse-of-dimensionality-free numerical approaches for solving
  hamilton–jacobi equations}, Applied Mathematics and Optimization, 83
  (2021), pp.~1--49.

\bibitem{YB_2012}
{\sc T.~T. Yusuf and F.~Benyah}, {\em Optimal control of vaccination and
  treatment for an {SIR} epidemiological model}, World journal of modelling and
  simulation, 8 (2012), pp.~194--204.

\end{thebibliography}

\end{document}